\documentclass{elsarticle}
\usepackage{amsmath,amsthm,amsfonts,latexsym,amssymb,fancybox,comment,graphicx,subfigure,epsfig,eepic,epic,epstopdf,array,subeqnarray}
\usepackage{calrsfs,mathrsfs,amsbsy}
\DeclareMathAlphabet{\pazocal}{OMS}{zplm}{m}{n}
\usepackage{colortbl,color}
\usepackage[colorlinks,linktocpage,linkcolor=blue]{hyperref}
\usepackage{hypernat,mathrsfs}

\def\emp#1{{\color{red}{\bf #1}}}
\usepackage[all]{xy}

\journal{}

\begin{document}
\newtheorem{theorem}{Theorem}[section]
\newtheorem{lemma}[theorem]{Lemma}
\newtheorem{proposition}[theorem]{Proposition}
\newtheorem{remark}[theorem]{Remark}
\newtheorem{corollary}[theorem]{Corollary}
\newtheorem{definition}[theorem]{Definition}
\newtheorem{example}[theorem]{Example}
\newtheorem{assumption}[theorem]{Assumption}

\newcommand{\thmref}[1]{Theorem~\ref{#1}}
\newcommand{\thmrefs}[2]{Theorems~\ref{#1} and ~\ref{#2}}
\newcommand{\secref}[1]{\S\ref{#1}}
\newcommand{\lemref}[1]{Lemma~\ref{#1}}
\newcommand{\lemrefs}[2]{Lemmas~\ref{#1} and ~\ref{#2}}
\newcommand{\propref}[1]{Proposition~\ref{#1}}
\newcommand{\rmkref}[1]{Remark~\ref{#1}}
\newcommand{\assref}[1]{Assumption~\ref{#1}}
\newcommand{\exrefs}[2]{Examples~\ref{#1} and ~\ref{#2}}
\newcommand{\exref}[1]{Example~\ref{#1}}
\newcommand{\figref}[1]{Fig.~\ref{#1}}
\newcommand{\figrefs}[2]{Figs.~\ref{#1} and ~\ref{#2}}
\newcommand{\tabref}[1]{Table~\ref{#1}}
\newcommand{\tabrefs}[2]{Tables~\ref{#1} and ~\ref{#2}}

\newcommand{\rb}[1]{\raisebox{1.5ex}[0pt]{#1}}
\newcommand{\Mid}[2]{\left.#1\right|_{#2}}
\def\<{\left\langle}
\def\>{\right\rangle}
\def\O{\Omega}
\def\G{\Gamma}
\def\p{\partial}
\def\hat{\widehat}
\def\hbx{\widehat{\mathbf x}}
\def\hx{\widehat x}
\def\hy{\widehat y}
\def\hQ{\widehat Q}
\def\ba{\mathbf a}
\def\mbb{\mathbf b}
\def\mbP{\mathbf P\textsc{}}
\def\mbQ{\mathbf Q}
\def\mbB{\mathbf B}
\def\bB{\mathbf B}
\def\mbe{\mathbf e}
\def\mbf{\mathbf f}
\def\bg{\mathbf g}
\def\dsum{\displaystyle\sum}
\def\dx{\operatorname{\,d\mathbf x}}
\def\ds{\operatorname{\,d{s}}}
\def\bn{\mathbf n}
\def\bp{\mathbf p}
\def\mcO{\mathcal O}
\def\bq{\mathbf q}
\def\br{\mathbf r}
\def\bs{\mathbf s}
\def\bt{\mathbf t}
\def\bu{\mathbf u}
\def\bug{\mathbf u_\bg}
\def\bv{\mathbf v}
\def\bw{\mathbf w}
\def\bx{\mathbf x}
\def\by{\mathbf y}
\def\bz{\mathbf z}
\def\bF{\mathbf F}
\def\bH{\mathbf H}
\def\bL{\mathbf L}
\def\bP{\mathbf P}
\def\bV{\mathbf V}
\def\bW{\mathbf W}
\def\ptl{\partial}
\def\bbP{\mathbb P}
\def\bbR{\mathbb R}
\def\Tau{\mathcal T}
\def\bPhi{{\boldsymbol \Phi}}
\def\bnu{{\boldsymbol \nu}}
\def\bphi{{\boldsymbol \phi}}
\def\bpsi{{\boldsymbol \psi}}
\def\bPi{{\boldsymbol \Pi}}
\def\Grad{\operatorname{grad}\,}
\def\Div{\operatorname{div}\,}
\def\Curl{\operatorname{curl}\,}
\def\grad{\nabla\,}
\def\div{\nabla\cdot}
\def\curl{\nabla\times}
\def\calL{\mathcal L}
\def\NC{\mathcal{NC}}
\def\bNC{\mathcal{\mathbf{NC}}}
\def\til{\widetilde}

\def\DSSY{DSSY}
\def\vDSSY{\mathbf{DSSY}}
\def\DSSYhz{{DSSY}^h_0}
\def\vDSSYhz{\mathbf{DSSY}^h_0}
\def\vcQhz{\pmb{\mathscr{Q}}_{1,0}^{c,h}} 
\def\vcQqhz{\pmb{\mathscr{Q}}_{2,0}^{c,h}} 
\def\dcQlhz{{\mathscr{Q}}_{1,0}^{dc,h}} 
\def\NChz{{\mathscr{P}}_{1,0}^{nc,h}}
\def\vNCh{\pmb{\mathscr{P}}_1^{nc,h}}
\def\vNChz{\pmb{\mathscr{P}}_{1,0}^{nc,h}}
\def\ttvNChz{\widetilde{\widetilde{\vNChz}}}
\def\msvNChz{\pmb{\mathscr{P}}_{1,0}^{me,2h}} 
\def\bp{\overline{p}_h}
\def\bq{\overline{q}_h}
\def\gcbs{{\mathcal C}^h}   
\def\gcb{{\mathbb C}_h}   
\def\gbbs{{\pmb{\mathscr B}}^h}   
\def\gbb{{\mathbb B}_h}   
\def\cdfs{{\pmb{\mathscr D}}}   
\def\gdfs{{\pmb{\mathscr D}}^h}   
\def\ttgdfs{\widetilde{\widetilde{\pmb{\mathscr D}^h}}}   
\def\Pch{{\mathcal P}_{c}^h}
\def\Phz{{\mathscr{P}}_0^{h}}
\def\PhRz{{\mathscr{P}}_{R,0}^{h}}
\def\PhBz{{\mathscr{P}}_{B,0}^{h}}
\def\Pcf{\widetilde{\Phz}}
\def\Pchp{{\mathcal P}_{c'}^h}
\def\Pcfp{{\mathcal P}_{cf'}^h}

\def\div{\nabla\cdot\,}
\newcommand{\beq}{\begin{equation}}
\newcommand{\eeq}{\end{equation}}
\newcommand{\beqs}{\begin{equation}}
\newcommand{\eeqs}{\end{equation}}
\newcommand{\beqx}{\begin{equation*}}
\newcommand{\eeqx}{\end{equation*}}
\newcommand{\bal}{\begin{aligned}}
\newcommand{\eal}{\end{aligned}}
\newcommand{\bes}{\begin{eqnarray*}}
\newcommand{\ees}{\end{eqnarray*}}
\def\supp{\operatorname{supp}}
\def\Span{\operatorname{Span}}
\newcommand{\wh}{\widehat}

\begin{frontmatter}

\title{Nonconforming finite element method applied to the driven cavity problem}
\tnotetext[label1]{The project is supported in part by
National Research Foundation of Korea (NRF--2014R1A2A1A11052429).}

\author[snucst]{Roktaek Lim}
\ead{rokt.lim@gmail.com}

\author[snucst,snumath]{Dongwoo Sheen}
\ead{dongwoosheen@gmail.com}


\address[nanyang]{School of Physical \& Mathematical Sciences, Nanyang Technological University,
180 Ang Mo Kio Ave 8, Singapore 569830}

\address[snumath]{Department of Mathematics and Interdisciplinary Program in
  Computational Science \& Technology, Seoul National University, Seoul 08826, Korea}

\begin{abstract}
A cheapest stable nonconforming finite element method is presented for
solving the incompressible flow in a square cavity without smoothing the corner singularities.
The stable cheapest nonconforming finite element pair based on $P_1\times P_0$ on rectangular meshes \cite{stab-cheapest} is employed with a minimal modification of the discontinuous Dirichlet data on the top boundary, where $\Pcf$ 
is the finite element space of piecewise constant pressures with the
globally one-dimensional checker-board pattern subspace eliminated. 
The proposed Stokes elements have the least number of degrees of freedom compared
to those of known stable Stokes elements. Three accuracy indications
for our elements are analyzed and numerically verified. Also, various
numerous computational results obtained by using our proposed element show excellent accuracy.
\end{abstract}


\end{frontmatter}

\section{Introduction} \label{sec:intro}
The lid driven square cavity has been one of the most popular benchmark
problems for new numerical methods for the incompressible Navier-Stokes
equations in terms of accuracy, numerical efficiency and so on. To refer only
few see
\cite{auteri2002numerical, botella1998benchmark, erturk2005numerical, erturk2009discussions},
for instance, and the references therein.
The presence of singularities at the upper corners of the cavity is the source
of numerical difficulties for solving the cavity flow problem. It is usually
erroneous to use high-order methods without handling the corner singularities due
to the Gibbs phenomenon. Many studies have been carried out to overcome this difficulty.
Barragy and Carey \cite{barragy1997} used a $p$-version finite element
formulation ($p\ge 6$) combined
with a strongly graded and refined mesh to handle the corner singularities.
Other studies change the boundary condition to
overcome this difficulty: see, for instance,
\cite{hna09, shen1991, sahin2003, glowinski-semicircular-06}, and the references therein.
The latter approach are coined as the so-called
regularized lid driven cavity problem. The constant boundary condition for velocity is replaced by
a function that vanishes at the upper corners of cavity \cite{hna09, shen1991}.
Botella and Peyret \cite{botella1998benchmark} solved a regularized cavity problem
by using a subtraction method of the leading terms from the asymptotic expansion
of the solution of the Navier-Stokes equations in the vicinity of the corners,
where the velocity is discontinuous.
Sahin and Owens \cite{sahin2003} inserted leaks across the heights of the
finite volumes at the corners between the lid and the vertical walls to handle
the corner singularities.
Many studies reported that in the critical Reynolds number range $[7000,8500]$
Hopf bifurcations occur for the lid driven square cavity problem
\cite{auteri2002numerical, erturk2005numerical, hna09, shen1991}.
Bruneau and Saad \cite{bruneau20062d} revisited the issue of bifurcation using
third--order time discretization schemes with the $5000\times 5000$ finite difference spatial
discretizations. They observed the first bifurcation occurs between $Re=8000$
and $Re=8050.$
Guermond and Minev \cite{guermond2012start} reported three--dimensional
benchmark solutions using a direction splitting method introduced in \cite{guermond2010new, guermond2011new}. 
They also provided two--dimensional solutions, which are
correct up to at least three digits, for $Re=1000$
using the uniform $5000\times 5000$ MAC stencil.
Instead of the square domain, Glowinski {\it et al.}
\cite{glowinski-semicircular-06} considered a semi-circular cavity-driven flow
with a special time-dependent regularization on the Dirichlet data at
the two corners: they observed Hopf bifurcations around $Re=6600$,
which is smaller than the case of square domain,
using an iso-parametric variant of the Bercovier-Pironneau element
\cite{bercovier-pironneau-79} introduced in \cite{hna09}.

The purpose of the current paper is to try to solve the lid driven square cavity problem
without any regularization at the corners, employing nonconforming finite 
element pairs whose degrees of freedom and implementation are as cheap as possible.
As the nonconforming elements use the values at the midpoints of edges as
DOFs, instead of those at the vertices, the discontinuity singularities at the corners
are naturally treated without any regularization. Our nonconforming finite
element pairs are based on the two stable nonconforming finite element pairs on
uniform square meshes \cite{stab-cheapest} introduced for the stationary
incompressible Stokes problem. The two pairs are briefly described as follows:
The first of them uses the $P_1$-nonconforming quadrilateral element \cite{parksheen-p1quad}
for the approximation of the velocity field, componentwise, while the pressure
is approximated by a subspace of the piecewise constant functions whose
dimension is two less than the number of squares in the mesh.
The second of them is a one-dimensional modification of the above finite
element pairs to both velocity and pressure spaces: the velocity space is
enriched by a globally one-dimensional DSSY(Douglas-Santos-Sheen-Ye)-type bubble function
\cite{dssy-nc-ell, cdy, jeon-nam-sheen-shim-nonpara} while the pressure space
is the subspace of the piecewise constant functions whose
dimension is one less than the number of squares in the mesh in order to
fulfill the mean-zero property.
The stability and optimal convergence results for these element
pairs applied to the stationary Stokes equations with the homogeneous Dirichlet
boundary condition can be found in \cite{stab-cheapest}.

In order to treat the inhomogeneous lid-driven Dirichlet boundary condition, 
we modified the above elements \cite{stab-cheapest} as follows. The boundary
condition on the interior of the top boundary is handled as usual, but the corner
boundary condition is specially treated at the two end elements on the top
by adding two local DSSY-type bubble
functions whose values at the midpoints of top boundary parts to be
$(1,0)^t$ and at the midpoint of the other boundary parts to be
$(0,0)^t$. Indeed, since nonconforming finite element methods can avoid 
vertex values degrees of freedom, the boundary values at the top left and right corners
are not required. Thus, one can solve the driven cavity problem
without any regularization of the boundary condition \eqref{cavityBD}. 

We note that the above modified finite elements have the smallest DOFs and are easiest
to implement among the finite element space pairs that fetch all non-spurious
piecewise constant pressure fields. Moreover, our finite element methods yield
nearly divergence free velocity fields. Indeed,
$\int_\O | \div {\mathbf v}_h | \, d\bx= \mcO (h^3)$,
which is good indication of numerical solver. The $\mcO (h^3)$ factor arises
from the inhomogeneous boundary data (the finite element pairs introduced in 
\cite{stab-cheapest} for the homogeneous boundary condition yield exactly
divergence free velocity approximation.) 
Another indication of superiority of our element is that our methods gives
substantially smaller volumetric flow rates across horizontal and vertical
line sections \cite{aydin} than other methods by
a factor of two. They are reported in \secref{sec:results}.

The plan of our presentation of this paper is as follows. In the next section 
the lid-driven cavity problem is briefly described.
With a brief review on the $P_1$-nonconforming quadrilateral element, a 
detailed description and implementation of our finite element methods are given
in \secref{sec:method}.
Three accuracy indications of our numerical solutions are analyzed in \secref{sec:acc}.
Some numerical results are presented \secref{sec:results} with comparison to
the results of other methods. The last section concludes our presentation.

\section{Problem formulation} \label{sec:prob}
Let $\O =(0,1)^2$ be the square cavity. Consider
the steady-state incompressible Navier-Stokes equations in dimensionless form:
\beq \label{nseqn}
\bal
- \nu \triangle \mathbf{u} +(\mathbf{u} \cdot \nabla)\mathbf{u} + \nabla p & = \mathbf{0}
\qquad \text{in} \; \O,\\
\div \mathbf{u} & = 0 \qquad \text{in} \; \O, \\
\eal
\eeq
with the Dirichlet boundary condition
\beq \label{bceqn}
\mathbf{u} = \mathbf{g} \quad \text{on} \; \G
\quad \text{with} \; \int_{\G} \bnu \cdot  \bg \ds = 0.
\eeq
Here, $\bu$ and $p$ denote the flow velocity and pressure, $\nu$ the fluid kinetic viscosity,
$\G$ the boundary of $\O$, and $\bnu$ the unit outward normal vector to $\O$.
Here, and in what follows, bold faces will denote the two-dimensional vectors,
functions, and function spaces.
For the driven cavity problem, suppose that the Dirichlet data is given by
\begin{equation} \label{cavityBD}
\mathbf{g}(x,y)=\left\{
\begin{aligned}
& (1,0)^t \; \text{ if } 0< x <1 \; \text{ and } \; y=1, \\
& \text{arbitrary } \; \text{ if } x=0, 1 \; \text{ and } \; y=1, \\
& \mathbf{0} \; \text{ elsewhere on } \partial \O.
\end{aligned}
\right.
\end{equation}
Notice that the regularity of the boundary value of velocity field:
$\mathbf{g} \in \bH^{\frac12-\epsilon}(\G)$  for
arbitrary $\epsilon > 0,$
which limits the regularity of the solution $\bu \in \bH^{1-\epsilon}(\O)$ at best. 
The possible highest regularity of the solutions is
\[
(\bu, p) \in \bW^{1,r}(\O) \times W^{0,r}(\O)/{\mathbb R} \quad\text{ for all }r \in (1,2).
\]
The Sobolev embedding theorem implies 
$
(\bu, p) \notin \bH^{1}(\O) \times L^{2}(\O)/{\mathbb R},
$
but
$
(\bu, p) \in \bH^{1-\epsilon}(\O) \times H^{-\epsilon}(\O)/{\mathbb R}
$
for arbitrary small $\epsilon > 0.$
See \cite{cai-wang-cavity} for more details of analysis in the case of driven cavity Stokes equations.

\section{A cheapest nonconforming finite element method} \label{sec:method}
In this section we will begin with a brief review on the $P_1$-nonconforming quadrilateral element
\cite{parksheen-p1quad, cpark-thesis, altmann-carstensen, altmann2012p}
Then the stable cheapest finite element pairs \cite{stab-cheapest} for the incompressible Stokes
equations with homogeneous boundary condition will be described. In the third
part of this section describes the treatment of nonhomogeneous boundary
condition for the lid-driven cavity problem. Especially the corner
singularities will be taken care of.

\subsection{The $P_1$-nonconforming quadrilateral element space}
In this paper, we consider the unit square domain $\O=(0,1)^2$ with uniform square meshes.
Let $(\Tau_h)_{0<h<1}$ be a family of partitions of $\O$ into $N_Q=N^2$
disjoint squares $Q_{jk}$ of size $h \times h$, $h=1/N,$ with
barycenter $\left( (j-\frac12)h,  (k-\frac12)h \right),$ for 
$j,k=1,\cdots,N.$ We assume that $N$ is an even integer.
By $N_{v}^{i}$ denote the number of interior vertices $V_{jk}=(jh,kh)$
in $\Tau_h$ so that $N_v^i = (N-1)^2.$
Set
\begin{equation*}
\begin{aligned}
\NChz=\{ & v \in L^{2}(\Omega) \; | \;
v|_{Q_{jk}} \in P_1(Q_{jk}) \; \forall Q_{jk} \in \Tau_h,\;
v \text{ is continuous at} \\
&\text{the mid point of each interior edge in } \Tau_h \text{ and } v
\text{ vanishes at} \\
&\text{the mid point of each boundary edge in } \Tau_h \},
\end{aligned}
\end{equation*}
The global basis functions of $\NChz$ can be defined
vertex-wise:
for each interior vertex $V_{jk}$ in $\Tau_h,$
define $\phi_{jk} \in \NChz$ such that it has value $1$ at the midpoint of each
interior edge whose end points contains the vertex $V_{jk}$ and
value $0$ at the midpoint of every other edge in $\Tau_h$. Then
the $P_1$-nonconforming quadrilateral element space
\cite{parksheen-p1quad, cpark-thesis} is given by
\[
\NChz=\left\{ v_h =  \sum_{j,k=1}^{N-1} \alpha_{jk}\phi_{jk}\,\mid
\alpha_{jk}  \in \mathbb R \quad \forall j,k \right\},\quad\dim(\NChz) = N_v^i.
\]

\subsection{The DSSY-type finite element space}
The $DSSY$ nonconforming element space on a reference domain $\hQ:=[-1,1]^2,$ 
with vertices $\hbx_1 = (1,0), \hbx_2 = (0,1), \hbx_3 = (-1,0), \hbx_4 = (0,-1),$
is defined by
$$
DSSY(\hQ)=\Span\{1,\hx,\hy, \theta_k(\hx)-\theta_k(\hy)\},
$$
where
$$
\theta_\ell(t) = \begin{cases}  t^2,\quad   &\ell  = 0, \\
 t^2-\frac{5}{3}t^4,\quad   &\ell  = 1, \\
 t^2-\frac{25}{6}t^4+\frac{7}{2}t^6 ,\quad   &\ell  = 2.
\end{cases}
$$
The reference DSSY basis functions  have the form
\begin{equation*}
\begin{aligned}
&\wh{\psi}_{\hbx_1}^{\text{DSSY}}(\hbx) =
\frac{1}{4}+\frac{1}{2}\hx-\frac{3}{8}(\theta_\ell(\hx)-\theta_\ell(\hy)),\\
&\wh{\psi}_{\hbx_2}^{\text{DSSY}}(\hbx) =
\frac{1}{4}+\frac{1}{2}\hy+\frac{3}{8}(\theta_\ell(\hx)-\theta_\ell(\hy)),\\
&\wh{\psi}_{\hbx_3}^{\text{DSSY}}(\hbx) =
\frac{1}{4}-\frac{1}{2}\hx-\frac{3}{8}(\theta_\ell(\hx)-\theta_\ell(\hy)),\\
&\wh{\psi}_{\hbx_4}^{\text{DSSY}}(\hbx) =
\frac{1}{4}-\frac{1}{2}\hy+\frac{3}{8}(\theta_\ell(\hx)-\theta_\ell(\hy)),
\end{aligned}
\end{equation*}
such that 
$\wh{\psi}_{\hbx_j}^{\text{DSSY}}(\hbx_k)=\delta_{jk},$ the Kronecker
delta. In what follows, we fix $\ell = 1.$

Let $F_Q:\hQ \rightarrow Q$ be a bijective affine transformation from the
reference domain onto a rectangle $Q$. Then $DSSY(Q)$ is defined by
\begin{eqnarray}
DSSY(Q) = \left\{ \hat{v}~\circ~F_Q^{-1}~\middle|~\hat{v}
\in DSSY(\hat{Q}) \right\}.
\end{eqnarray}
Then the DSSY-type finite element space
\cite{cdssy,cdy,dssy-nc-ell,jeon-nam-sheen-shim-nonpara} is defined by
\begin{eqnarray*}
&&{DSSY^h_0}=\{v\in L^2(\O)~|~ v|_Q\in
                 DSSY(Q) ~ \forall Q\in \Tau_h; \\
&&\qquad\qquad  \qquad  v \text{ is continuous at the midpoint of each interior edge }  \\
&&\qquad\qquad  \qquad ~~ \text{ and vanishes at the midpoint of each boundary
                   edge in }
 \Tau_h \}.
\end{eqnarray*}
\begin{remark} For $\ell=0,$ the DSSY-type nonconforming element, or the
velocity components in the CDY(Cai-Douglas-Ye) Stokes element, is identical to the
  rotated $Q_1$ element of Rannacher and Turek \cite{rann}. The difference
  between and the DSSY-type nonconforming elements (with $\ell=2,3$) and the
  rotated $Q_1$ element is that the former satisfies the mean value property
  $\frac1{|e|}\int_e v\ds = v(m_e)$ on each edge $e$ in $\Tau_h,$ where $m_e$
  denotes the midpoint of $e.$ See \cite{jeon-nam-sheen-shim-nonpara} for more details.
\end{remark}

\subsection{The stable cheapest finite element pairs: homogeneous Dirichlet boundary case}
Set 
\[
\Phz = \left\{p_h=\sum_{j,k=1}^N \gamma_{jk}\chi_{Q_{jk}}\,\mid\, \gamma_{jk} \in
\mathbb R;\, \int_\O p_h\,\dx = 0  \right\} \subset L^2_0(\O),
\]
where $\chi_{Q_{jk}}$ denotes the usual characteristic function.
Denote by $\Pcf$ the subspace of $\Phz$ by removing the globally
one-dimensional global checkerboard pattern from $\Phz.$
One way of forming the basis for the $(N_Q-2)$--dimensional space $\Pcf$
can be described as follows. 
Let $\O =\O^R \cup \O^B$ be a decomposition of $\O$ into the disjoint
unions of red and black rectangles 
$\O^R = \cup_{Q\in \Tau_h^R } Q$ and 
$\O^B = \cup_{Q\in \Tau_h^B } Q,$ where
\begin{eqnarray*}
\Tau_h^R &=& \{Q_{jk}\in \Tau_h\,\mid\quad j+k\text{ is an even integer}\},\\
\Tau_h^B &=& \{Q_{jk}\in \Tau_h\,\mid\quad j+k\text{ is an odd integer}\}. 
\end{eqnarray*}

We are now in a position to form the two $\left(\frac{N_Q}{2}-1\right)$--dimensional
subspaces of $\Phz$ as follows:
\begin{eqnarray*}
\PhRz &=& \left\{p_h=\sum_{j,k=1}^N \gamma_{jk}\chi_{Q_{jk} \cap \O^R }\,\mid\, \gamma_{jk} \in
\mathbb R;\, \int_\O p_h\,\dx = 0  \right\} \subset L^2_0(\O),\\
\PhBz &=& \left\{p_h=\sum_{j,k=1}^N \gamma_{jk}\chi_{Q_{jk} \cap \O^B }\,\mid\, \gamma_{jk} \in
\mathbb R;\, \int_\O p_h\,\dx = 0  \right\} \subset L^2_0(\O).
\end{eqnarray*}
Then it turns out that $\Pcf = \PhRz\oplus \PhBz,$ from which the basis
functions
for $\Pcf$ is built in a standard way by taking the union of the basis functions of 
$\PhRz$ and $\PhBz.$
Henceforth the first pair of stable cheapest finite element pair for the
incompressible Stokes flows is given as
\begin{eqnarray}\label{eq:1st-fempair}
\vNChz \times \Pcf\quad\text{with dimension }
2N^i_v+N_Q-2.
\end{eqnarray}
A second pair of stable cheapest finite element pair is obtained by enriching
the velocity space by a globally one-dimensional
Assume that $N$ is an even integer.
Denote by $\Tau^h_M$ the macro mesh such that each macro rectangle $Q_{JK}^M$ 
consists of $2\times 2$ rectangles $Q_{jk},Q_{j,k+1}, Q_{j+1,k}, Q_{j+1,k+1},$ with $(J,K) =(j,k).$ 
from $\Tau^h$  with $J,K = 1,3,\cdots, N-1.$
For each macro-element $Q_{JK}^M \in \Tau^h_M$, define $\bpsi_{Q_{JK}^M} \in \mathbf{DSSY}_0^h$
such that
\bes
\supp(\bpsi_{Q_{JK}^M} ) \subset \overline{Q}_{JK}^M,
\ees
and the integral averages over the edges in $\Tau_h$ vanish except
\bes
\oint_{\p Q_{j,k} \cap \p Q_{j+1,k}} \bpsi_{Q_{JK}^M}~ds = \bnu,\qquad
\oint_{\p Q_{j,k+1} \cap \p Q_{j+1,k+1}} \bpsi_{Q_{JK}^M}~ds = -\bnu.
\ees
where $\bnu$ denotes the unit outward normal vector to $Q_{j\ell}$ on the edge
$\p Q_{j\ell} \cap \p Q_{j+1\ell}$, $\ell=k,k+1$. 
Introduce the following vector space  of macro bubble functions:
$\gbbs=\Span\left\{ \dsum_{Q_{JK}^M\in\Tau^M}\bpsi_{Q_{JK}^M}\right\} 
$
which is a one-dimensional subspace of $\mathbf{DSSY}^h_0.$
Then $\vNChz$ in enriched by adding $\gbbs,$ denoted by
$\ttvNChz=\vNChz\oplus \gbbs$,
and hence, the second stable Stokes finite element pair is defined as follows:
\begin{eqnarray}\label{eq:2nd-fempair}
\ttvNChz \times \Phz\quad\text{with dimension }
2N^i_v+N_Q.
\end{eqnarray}

The stability and optimal convergence properties of the two pairs of Stokes elements 
\eqref{eq:1st-fempair} and \eqref{eq:2nd-fempair}
are shown for the stationary Stokes equations in \cite{stab-cheapest}.

Comparing several other stable quadrilateral finite element pairs satisfying
the inf-sup condition \cite{crouzeix-raviart, taylor-hood, rann},
the nonconforming element pairs \eqref{eq:1st-fempair} and \eqref{eq:2nd-fempair}
 have the lowest degrees of freedom.
\tabref{fedofcomp} illustrates the degrees of freedom for different pairs,
whose notations will be used throughout the paper.

\begin{table}
\caption{Number of degrees of freedom for different pairs
(velocity/pressure)}
\centering
\begin{tabular}{ccccc}
\hline
$N$ & $\vcQqhz\times \dcQlhz$ & $\vcQqhz\times \Phz$ & $\text{\bf rot} \mbQ_1
\times \Phz$ & $\vNChz \times \Pcf$
  \\ \hline
$2^{4}$ & 2178/289 & 2178/255 & 1088/255 & 450/254 \\
$2^{5}$ & 8450/1089 & 8450/1023 & 4224/1023 & 1922/1022 \\
$2^{6}$ & 33282/4225 & 33282/4095 & 16640/4095 & 7938/4094 \\
$2^{7}$ & 132098/16641 & 132098/16383 & 66048/16383 & 32256/16382 \\
\hline
\end{tabular}\caption{ $\vcQqhz\times \dcQlhz$, $\vcQqhz\times \Phz$, $\text{\bf rot} \mbQ_1
\times \Phz$, and $\vNChz \times \Pcf$ stand for the two Taylor--Hood elements
of type
$Q_2\times Q_1$ and $Q_2\times P_1$, the Rannacher--Turek  nonconforming quadrilateral rotated $Q_1\times
P_0$ element, and the nonconforming quadrilateral $P_1\times P_0$ element \cite{stab-cheapest}, respectively.}
\label{fedofcomp}
\end{table}
It is shown that both nonconforming finite element spaces
\eqref{eq:1st-fempair} and \eqref{eq:2nd-fempair} give exactly identical
solutions for velocity fields but slight different pressure solutions whose
differences in $L^2(\O)$-norm are of order $\mathcal O(h)$, and thus both
velocity and pressure are approximated with optimal convergence for Stokes flows.
See \cite[\S
  4]{stab-cheapest} for details. Due to this observation, we concentrate on the
finite element pair \eqref{eq:1st-fempair} for approximating the cavity flow.

\subsection{Treatment of nonhomogeneous boundary condition}
In order to deal with the
Dirichlet boundary values of cavity flows, the open boundary part (top
boundary) is modified with two additional DSSY-type elements located at the
two top corners,  ${Q_{1N}}$ and ${Q_{NN}}.$ Notice that for $j=1,\cdots, N-1,$
$\phi_{j,N}$ 
has value 1 at the midpoints $((j-\frac12)h, 1)$ and $((j+\frac12)h, 1)$ and 0
at the other midpoints on the top boundary, one sees that
\[
\binom{\frac12}{0} \sum_j^{N-1}  \phi_{j,N}
\] 
assigns the vector value $(1,0)^t$ at the midpoints
$((1+\frac12)h, 1),\cdots,((N-1-\frac12)h, 1)$ and 
$(\frac12,0)^t$ at the midpoints $(\frac12 h, 1)$ and $(1-\frac12 h, 1),$
respectively. Denote the DSSY basis functions whose supports are the top
corner elements as follows:
\[
\psi_{2,TL} = \left(\wh{\phi}_{\hbx_2}^{\text{DSSY}}\circ
F_{Q_{1N}}^{-1}\right)\quad\text{and }
\psi_{2,TR}= \left(\wh{\phi}_{\hbx_2}^{\text{DSSY}}\circ
F_{Q_{NN}}^{-1}\right),
\]
both of which have values 1 at the top midpoints $(\frac{1}{2} h, 1)$ and
$(1-\frac{1}{2} h, 1),$ respectively, and 0 at the other midpoints of the two elements.
Summarizing the above, the approximate nonconforming finite element solution
with the Dirichlet boundary data for the lid-driven cavity flow
is approximated by $\bu_h$ of the form
\begin{eqnarray} \label{appu}
\bu_h  &=& \bu_{0,h} + \bu_{b,h},\qquad\text{where} \\
\bu_{0,h} &=& \sum_{j,k=1}^{N-1} \binom{\xi_{jk}}{\eta_{jk}} \phi_{jk}, \quad
\bu_{b,h} = \binom{ \frac12  }{0} \left[ \sum_{j=1}^{N-1} \phi_{j,N} +\psi_{2,TL} 
+\psi_{2,TR} \right].\nonumber
\end{eqnarray}

We are now in a position to define a discrete weak formulation of \eqref{nseqn}
to find $(\bu_{0,h}, p_h) \in  \vNChz \times \Pcf$ such that
\begin{subeqnarray} \label{dweaknseqn}
a_h( \bu_{0,h}, \bv_h)
&+&c_h(\bu_{0,h} ; \bu_{0,h}, \bv_h)
+c_h(\bu_{0,h} ; \bu_{b,h}, \bv_h)
+c_h(\bu_{b,h} ; \bu_{0,h}, \bv_h)
\nonumber \\
+b_h(\bv_h, p_h )&=&-a_h( \bu_{b,h}, \bv_h) -c_h(\bu_{b,h} ; \bu_{b,h}, \bv_h) \qquad
\forall \bv_h \in [\NChz]^{2}, \slabel{dweaknseqn-a} \\
b_h(\bu_{0,h}, q_h ) & = & -b_h(\bu_{b,h}, q_h ) \qquad \forall q_h
\in \Pcf,
\slabel{dweaknseqn-b}
\end{subeqnarray}
where
\begin{gather*}
a_h(\bu_h,\bv_h) = \nu \sum_{Q\in\Tau_h} \int_Q \nabla \bu_h : \nabla \bv_h\dx,\qquad
b_h(\bv_h,q_h) = - \sum_{Q\in\Tau_h} \int_Q (\nabla \cdot \bv_h)q_h\dx,\\
c_h(\bw_h;\bu_h,\bv_h) =  \sum_{Q\in\Tau_h} \int_Q (\bw_h \cdot \nabla) \bu_h \cdot \bv_h \dx.
\end{gather*}
The nonlinear equations \eqref{dweaknseqn} can be approximated by the Picard iteration method
\cite{cuvelier1986finite, elman2014finite, karakashian1982galerkin}.
With an initial guess 
$(\bu^{(0)}_{0,h},p^{(0)}_h) \in \vNChz\times \Pcf$, define
the Picard iterates
$(\bu^{(k)}_{0,h}, p^{(k)}_h)\in  \vNChz\times \Pcf$ 
for $k=1,2,\cdots,$  solving the following Oseen problem: 
\begin{subeqnarray} \label{dweaknseqn-picard}
&&a_h(\bu_{0,h}^{(k)}, \bv_h)
+c_h(\bu_{0,h}^{(k-1)} + \bu_{b,h} ; \bu_{0,h}^{(k)}, \bv_h)
 +b_h(\bv_h, p_h^{(k)} )
\nonumber \\
&&\quad=-a_h( \bu_{b,h}, \bv_h) -c_h(\bu_{0,h}^{(k-1)} + \bu_{b,h} ; \bu_{b,h}, \bv_h)
 \quad \forall \bv_h \in \vNChz, \\
&&b_h(\bu_{0,h}^{(k)}, q_h )  =  -b_h(\bu_{b,h}, q_h ) \quad \forall q_h \in \Pcf.
\end{subeqnarray}
The Picard iterates
$(\bu^{(k)}_{0,h}, p^{(k)}_h)_{k\ge 1}$ are shown to
converge at a linear order to the solution $(\bu_{0,h}, p_h)$ of
\eqref{dweaknseqn} in \cite{karakashian1982galerkin}.
One may of course use the Newton iterates which converge quadratically with
sufficiently close initial guesses to the exact solution as described in
\cite{cuvelier1986finite, elman2014finite, karakashian1982galerkin}.

\section{Accuracy of solutions} \label{sec:acc}
Previous studies validated their numerical solutions by comparing
their numerical results with benchmark solutions in the literature,
for example, \cite{barragy1997} and \cite{ghia}.
According to Erturk {\it et al.} \cite{erturk2005numerical}, there are many
different numerical procedures for the lid-driven cavity
flow problem which yield very similar numerical results in the case of
$\text{Re} \leq 1000$,
however, their numerical solutions start to deviate from each other as 
the Reynolds number increases.

Hence, in order to claim some sort of superiority of our nonconforming method
over the other existing methods, we will not only compare our numerical results
with those in the literature, but also show some other assessments for the
accuracy of the numerical solution.

\subsection{Volumetric flow rate}
Aydin and Fenner \cite{aydin} suggested a measurement of the accuracy of
numerical solutions. They computed the net
volumetric flow rate, $Q$, passing through a vertical line and a horizontal
line to check the continuity of the fluid. Denote $\bu=(u,v)$, and
let $Q_{u,c}$ and $Q_{v,c}$ be the volumetric flow rate passing through a
vertical line $x=c$ and a horizontal line $y=c$, respectively.
The volumetric flow rate values, $Q_{u,c}$ and $Q_{v,c}$ can be computed by
\beq\label{eq:volumetric}
Q_{u,c}= \left| \int_0^{1} u(c,y) \; dy \right|,\qquad
Q_{v,c}= \left| \int_0^{1} v(x,c) \; dx \right|.
\eeq

\subsection{Compatibility condition for the stream function $\psi$}
We can also use the compatibility condition for the stream function $\psi$
for the assessment of the accuracy of numerical solutions. 
Using the expressions of the vorticity $\omega$ as the two-dimensional curl of the velocity:
$
\omega=\curl\bu,
$
and the velocity field as the two-dimensional curl of the stream function
$
\bu = \curl \psi,
$
one has the Neumann boundary value problem for $\psi$ as follows:
\beq
-\triangle \psi = \omega\quad \text{ in } \O, 
\eeq
with 
\[
\frac{\p \psi}{\p n} =\begin{cases} -u & \text{ for } y=0, 0<x<1, \\
u & \text{ for } y=1, 0<x<1, \\
v & \text{ for } x=0, 0<y<1, \\
-v & \text{ for } x=1, 0<y<1.
\end{cases}
\]
A compatibility condition, combined with \eqref{bceqn}, yields
\beq \label{vorcondition}
\int_{\O} \omega  \; \dx = -\int_{\ptl \O} \frac{\p\psi}{\p n} \; \ds =
-\int_{0}^{1} 1 \; \dx = -1. 
\eeq
One can compute $\int_{\O} \omega \; \dx$ by using
the numerical solution $\bu_h$, and compare to check
the accuracy of the numerical approximation.

\subsection{Incompressibility condition}
Since the pointwise incompressible condition
$
\div \bu=0
$
should hold pointwise, the smallness of
\beq \label{qdiv}
\max_{Q_{jk} \in \Tau_h} \left|
\int_{Q_{jk}} \div \bu_h \dx
\right|
\eeq 
is a good indicator to check numerical accuracy. This implies that 
\eqref{qdiv} of the numerical solution
$\bu_h$ should be close to zero.

Invoking \eqref{appu}, and observing that
\begin{subeqnarray}\label{eq:div0-dssy}
\sum_{Q\in\Tau_h}\int_Q \div \binom{\psi_{2,TL}}{0}  \dx&=&
\int_{\p Q_{1N}} \bnu\cdot \binom{\psi_{2,TL}}{0} \ds=0,\\
\sum_{Q\in\Tau_h}\int_Q  \div \binom{\psi_{2,TR}}{0}  \dx &=&
\int_{\p Q_{NN}} \bnu\cdot \binom{\psi_{2,TR}}{0} \ds=0,
\end{subeqnarray}
one sees the following simplification:
\begin{eqnarray}\label{eq:div-red}
\sum_{Q\in\Tau_h}\int_Q   \div  \bu_h \dx =
\sum_{j,k=1}^{N-1} \int_{\O}  \xi_{jk} \frac{\p \phi_{jk}}{\p x}
+ \eta_{jk} \frac{\p \phi_{jk}}{\p y} \dx
+ \frac12 \sum_{j=1}^{N-1} \int_\O \frac{\p\phi_{j,N}}{\p x}\dx.
\end{eqnarray}
Recall that $\phi_{jk}$ is piecewise linear, and hence its derivative is
constant on each $Q_{\ell m}$: indeed,
\begin{eqnarray}\label{eq:div-phi}
\nabla \phi_{jk} =\begin{cases}
\frac{1}{h}\left(1,1\right)^t & \text{on} \;
Q_{jk}, \\
\frac{1}{h}\left(-1,1\right)^t & \text{on} \;
Q_{j+1,k}, \\
\frac{1}{h}\left(-1,-1\right)^t & \text{on} \;
Q_{j+1,k+1}, \\
\frac{1}{h}\left(1,-1\right)^t & \text{on} \;
Q_{j,k+1}, \\
\end{cases}
\end{eqnarray}
for $j=1,\cdots,N-1,k=1,\cdots,N.$ 
Set $q_h=\sum_{j,k=1}^N \zeta_{jk}\chi_{Q_{jk}} \in \Pcf$
for a general piecewise constant element.
By exploiting $|Q_{\ell m}|=h^2,$ 
from \eqref{dweaknseqn-b}, \eqref{eq:div0-dssy}, and \eqref{eq:div-phi} it then follows that
\begin{eqnarray*}
\int_{\O} (\div \bu_h)q_h \dx = 
h \sum_{j,k=1}^{N-1} \left[\xi_{jk}(\zeta_{jk}-\zeta_{j+1,k+1}) + \eta_{jk}
(\zeta_{j+1,k}-\zeta_{j,k+1})\right] + \frac{h}2 \sum_{j=1}^{N-1} \zeta_{jN}
\end{eqnarray*}
As a basis for $\Pcf$, choose the
union of the basis functions of $\PhRz$ and $\PhBz.$ 
For each $j,k=1,3,5,\cdots, N-1,$
with 
$q_h = \chi_{Q_{jk}} - \chi_{Q_{N,N}} \in \PhRz$
and
$q_h = \chi_{Q_{j+1,k+1}} - \chi_{Q_{N,N}} \in \PhRz$
one sees from \eqref{dweaknseqn-b}, \eqref{eq:div0-dssy} that
\begin{eqnarray} \label{alphabeta-1}
\int_{Q_{jk}} \div \bu_h\dx= \int_{Q_{j+1,k+1}} \div \bu_h\dx= \int_{Q_{NN}} \div \bu_h\dx.
\end{eqnarray}
Similarly, for $j,k=1,3,5,\cdots, N-1,$
with 
$q_h = \chi_{Q_{j+1,k}} - \chi_{Q_{N-1,N}} \in \PhBz$
and
$q_h = \chi_{Q_{j,k+1}} - \chi_{Q_{N-1,N}} \in \PhBz$
one concludes from \eqref{dweaknseqn-b}, \eqref{eq:div0-dssy} that
\begin{eqnarray} \label{alphabeta-2}
\int_{Q_{j+1,k}} \div \bu_h\dx= \int_{Q_{j,k+1}} \div \bu_h\dx= \int_{Q_{N-1,N}} \div \bu_h\dx.
\end{eqnarray}
%
Setting $\gamma_1=\int_{Q_{NN}} \div \bu_h\dx$ and  $\gamma_2=\int_{Q_{N-1,N}}
\div \bu_h\dx,$  one obtains from
\eqref{alphabeta-1}  and \eqref{alphabeta-2} that
\beq
\sum_{Q\in\Tau_h^R}\int_Q \div \bu_h \dx = \frac{N_Q}{2} \gamma_1, \quad
\text{and} \;
\sum_{Q\in\Tau_h^B}\int_Q \div \bu_h \dx = \frac{N_Q}{2} \gamma_2.
\eeq
Consequently,
\[
\sum_{Q\in\Tau_h}\int_Q \div \bu_h\dx =\sum_{Q\in\Tau_h^R}\int_Q \div \bu_h\dx
+ \sum_{Q\in\Tau_h^B}\int_Q \div \bu_h\dx = \frac{N_Q}{2}(\gamma_1+\gamma_2).
\]
However, using the Divergence Theorem piecewise for each $Q\in\Tau_h$, we have
\[
\sum_{Q\in\Tau_h}\int_Q \div \bu_h\dx =\sum_{Q\in\Tau_h}\int_{\p Q} \bnu\cdot
\bu_h\ds =
\sum_{j=1}^N\int_{\p Q_{jN} \cap \{y=1\}} \bnu\cdot \bu_h\ds = 0.
\]
Hence, $\frac{N_Q}{2}(\gamma_1+\gamma_2)=0$ and therefore
\begin{subeqnarray}\label{eq:div-red-black}
\int_{Q_{jk}} \div \bu_h \dx &=& \gamma_1\quad\forall Q_{jk} \in \Tau_h^R,\\
\int_{Q_{jk}} \div \bu_h \dx &=& -\gamma_1\quad\forall Q_{jk} \in \Tau_h^B.
\end{subeqnarray}
In order to compute $\gamma_1$ exactly, we sum $\int_{Q_{jk}} \div\bu_h\dx$
over all the red-type rectangles $Q_{jk}\in \Tau_h^R$ invoking the form of 
$\bu_h$ given in \eqref{appu}. 
Observing that for each interior vertex $V_{jk}$ there are two rectangles in
$\Tau_h^R$ which share only the vertex: these two rectangles can be either the
pair $(Q_{jk}, Q_{j+1,k+1})$ or the pair
$(Q_{j+1,k},Q_{j,k+1}).$ Then the integrals $\int \div
{\xi_{jk}\choose\eta_{jk}} \phi_{jk}\dx$ over those pairs
cancel each other due to the Divergence Theorem or direct integrations.
Recalling that the DSSY-type basis function parts on the two top corners do
not contribute anything for the integration of divergence, we see that
\begin{eqnarray}
\sum_{Q_{jk}\in \Tau_h^R}\int_{Q_{jk}} \div\bu_h\dx  
&= &
\sum_{Q_{jk}\in \Tau_h^R}\int_{Q_{jk}} \div (\bu_{0,h} + \bu_{b,h})\dx\nonumber \\
&= &
\sum_{Q_{jk}\in \Tau_h^R}\int_{Q_{jk}} \div \bu_{b,h} \dx\nonumber \\
&= &
\sum_{Q_{jN}\in \Tau_h^R}\int_{Q_{jN}} \div  \binom{ \frac12  }{0} \sum_{j=1}^{N-1} \phi_{j,N} \dx
\nonumber \\
&= &\sum_{j=1}^{N/2} \int_{Q_{2j,N}} \div  \binom{ \frac12  }{0} (\phi_{2j-1,N} + \phi_{2j,N})
\dx \nonumber \\
&= &\sum_{j=1}^{N/2} \int_{\p Q_{2j,N}} \bnu\cdot \binom{ \frac12  }{0}
(\phi_{2j-1,N} + \phi_{2j,N}) \ds\nonumber \\
&= & h \sum_{j=1}^{N-1} \left( - \frac12 \right )^{j}= -\frac{h}{2}.
\label{eq:div-red-final}
\end{eqnarray}
A combination of \eqref{eq:div-red-black} and \eqref{eq:div-red-final} shows that
\beq \label{qdiv-final}
\left| \int_{Q_{jk}} \div \bu_h \dx\right| = h^3 \quad\forall Q_{jk}\in \Tau_h.
\eeq

We summarize the above results as in the following theorem:
\begin{theorem}\label{thm:div-q}
Let $\bu_h$ be in the form \eqref{appu} and 
$(\bu_{0,h}, p_h)\in  \vNChz\times \Pcf$ fulfills \eqref{dweaknseqn-b}. Then
\eqref{qdiv-final} holds. Moreover, the signature of the integral over
$Q_{jk}$'s are alternating.
\end{theorem}

\section{Numerical results} \label{sec:results}
We have computed the steady state solutions of lid driven cavity flow from
$\text{Re}=100$ to $\text{Re}=5000$ by using the Picard iteration method 
with the termination condition:
\beq
\left\|
\begin{pmatrix}
\mbf - \nu A \bu^{(k)} - N\bu^{(k)} - B^{T} \bp^{(k)} \\
\bg - B \bu^{(k)}
\end{pmatrix} \right\|
\leq
10^{-10}
\begin{pmatrix}
\mbf \\ \bg
\end{pmatrix}
\eeq
where
$A, B,$ and $N$ denotes the matrices for the discrete Laplacian,  divergence, and
convection, respectively.
First, notice that our proposed nonconforming finite element method for
finding $(\bu_{0,h}, p_h) \in \vNChz \times \Pcf$ fulfilling \eqref{dweaknseqn}
does not modify the discontinuities at
the top corners. However, the usual conforming finite element methods require
suitable modifications. For instance, for the $\vcQqhz \times \dcQlhz$ element method,
the two popular cavity boundary conditions are used:
the watertight cavity boundary condition
\beq \label{watertight}
\bg=\left\{
\bal
& (1,0)^t, \; \text{ if } 0 < x < 1 \; \text{ and } \; y=1, \\
& \mathbf{0}, \; \text{ elsewhere on } \partial \O,
\eal
\right.
\eeq
and the leaky cavity boundary condition
\beq \label{leaky}
\bg=\left\{
\bal
& (1,0)^t, \; \text{ if } 0 \leq x \leq 1 \; \text{ and } \; y=1, \\
& \mathbf{0}, \; \text{ elsewhere on } \partial \O,
\eal
\right.
\eeq
respectively.
Notice that both conditions \eqref{watertight} and \eqref{leaky} satisfy \eqref{vorcondition}.
Our own FORTRAN and MATLAB codes were developed to implement
the $\vNChz \times \Pcf$ element method while
the IFISS S/W\cite{ifiss} was used to implement the $\vcQqhz \times \dcQlhz$ element method.
\begin{table}
\caption{Values used to plot the contours of the stream function and the vorticity}
\label{plotvalue}
\centering
\begin{tabular}{ll}
\hline
Contours & Values \\ \hline
Stream function &
-0.1175, -0.1150, -0.11, -0.1, -0.09, -0.07, -0.05, -0.03, -0.01,\\
&
-1.0E-04, -1.0E-05, -1.0E-07, -1.0E-10, 1.0E-08, 1.0E-07,\\
&
1.0E-06, 1.0E-05, 5.0E-05, 1.0E-04, 2.5E-04, 5.0E-04, \\
&
1.0E-03, 1.5E-03, 3.0E-03 \\
Vorticity &
-5.0, -4.0, -3.0, -2.0, -1.0, -0.5, 0.0, 0.5, 1.0, 2.0, 3.0, 4.0, 5.0\\
\hline
\end{tabular}
\end{table}

\begin{figure}[htb!]\label{fig:velocity-profile}
\caption{Comparison of the vertical components of the $u$-velocity
along the segment $x \in [0,1]$, $y=1/2$ and 
the horizontal components of the $v$-velocity
along the segment $x \in [0,1]$, $y=1/2$ 
with $Re=1000$}
\centering
\includegraphics[width=0.49\textwidth]{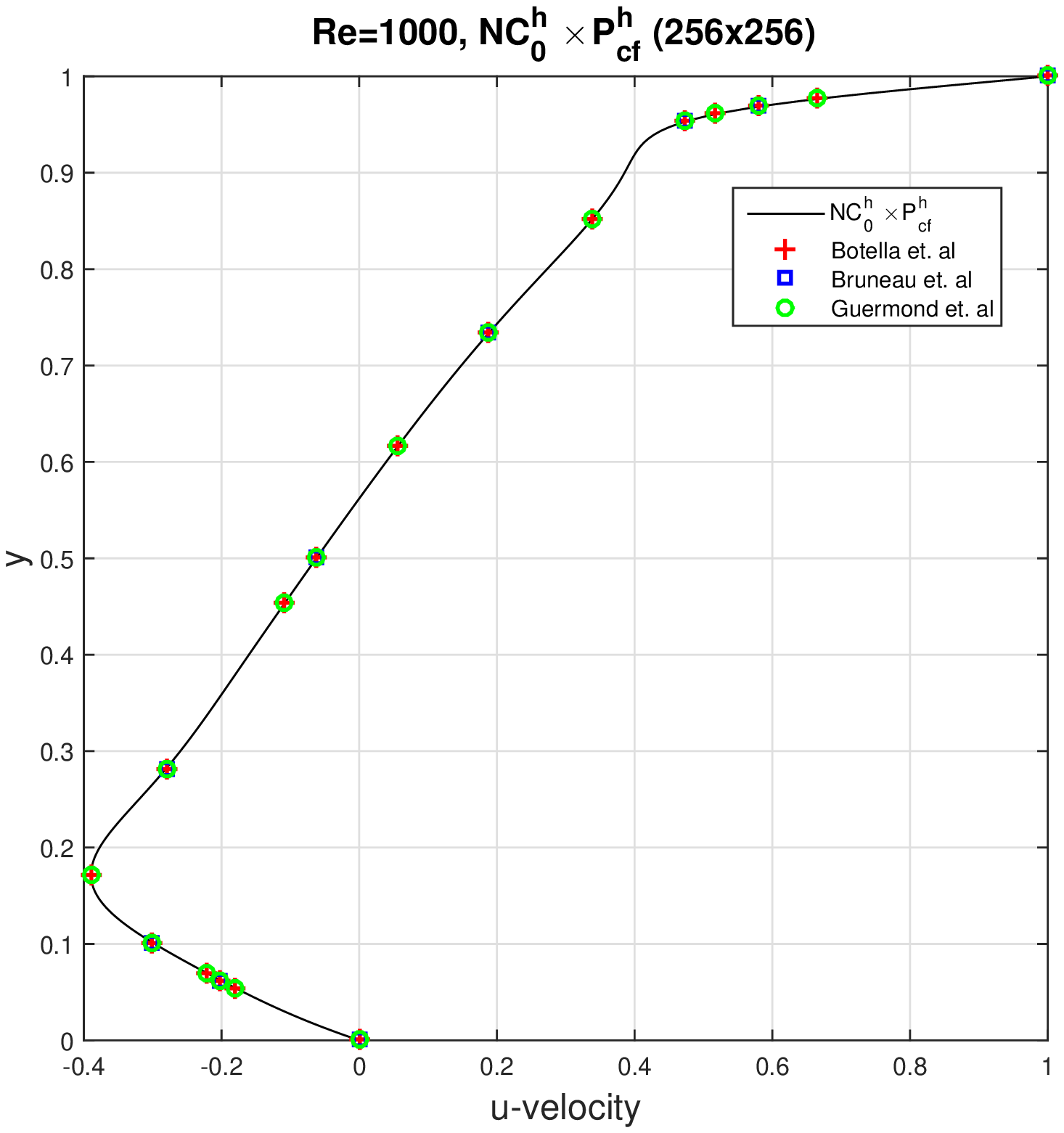}
\includegraphics[width=0.49\textwidth]{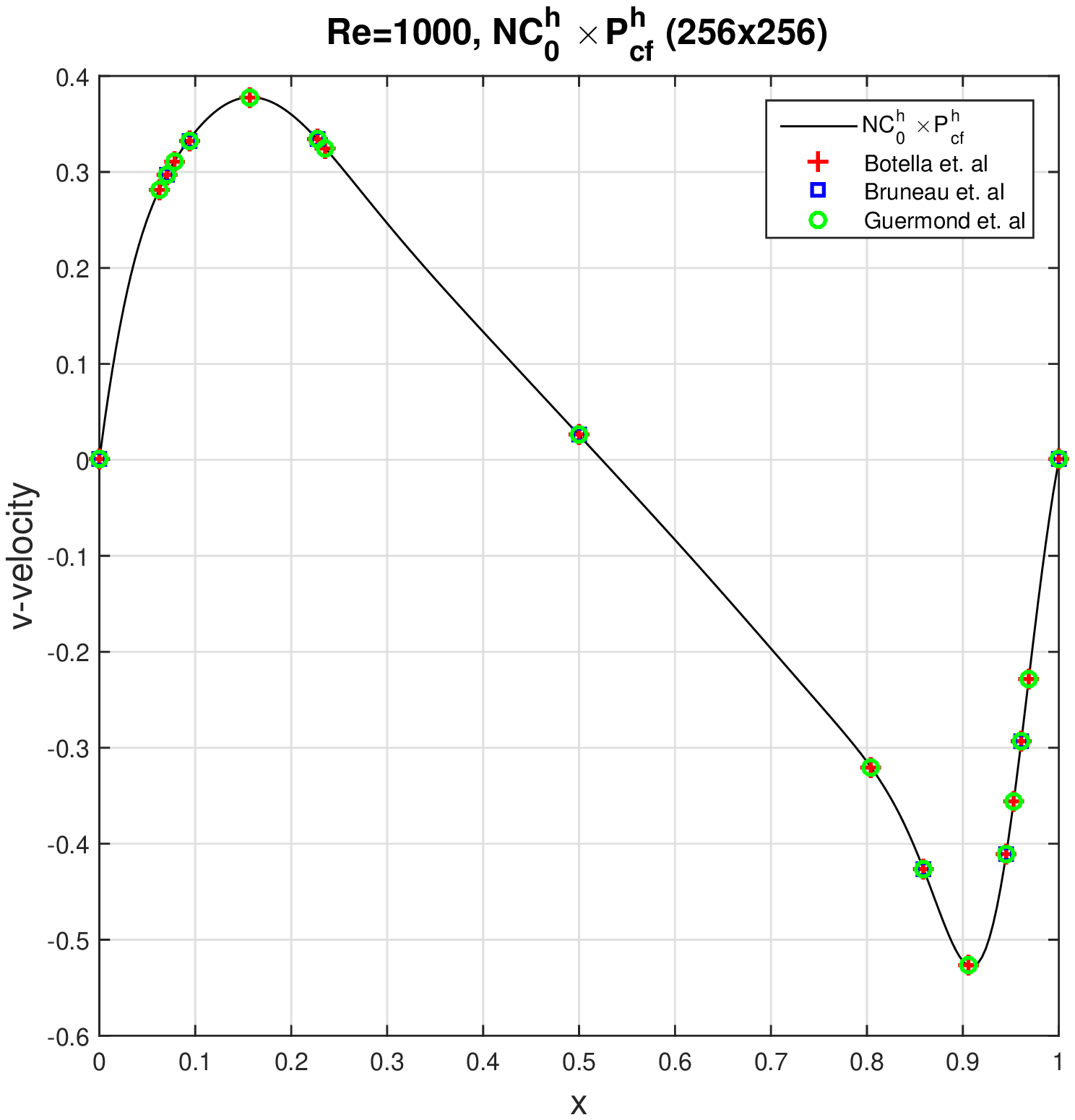}
\end{figure}

For the case of $Re=1000$, we present in Fig. \ref{fig:velocity-profile} the $u$-velocity profiles
along the line $x=0.5$ and the $v$-velocity profiles along the line $y=0.5$ computed
by using the $\vNChz \times \Pcf$ element with the boundary condition \eqref{cavityBD}
and compare our results with those by Botella and Peyret \cite{botella1998benchmark}, 
by Bruneau and Saad \cite{bruneau20062d}, and by Guermond and Minev \cite{guermond2012start}.
In each case, our velocity profiles show a good agreement with the
reference solutions.
Recall that the solutions obtained
Botella and Peyret used a spectral method on $160\times 160$ spectral nodes, 
and Bruneau and Saad used a second--order finite difference schemes on
the uniform $1024\times 1024$ nodes, while 
Guermond and Minev used a massively parallel computation combining
their new direction splitting algorithm and the MAC central finite difference scheme
on $5000\times 5000$ nodes.

\begin{figure}
\caption{Profiles of $u$-velocity along the line $x=0.5$ computed by using the stable
$\vNChz \times \Pcf$ with unregularized boundary condition}
\label{uprofile}
\centering
\includegraphics[width=0.49\textwidth]{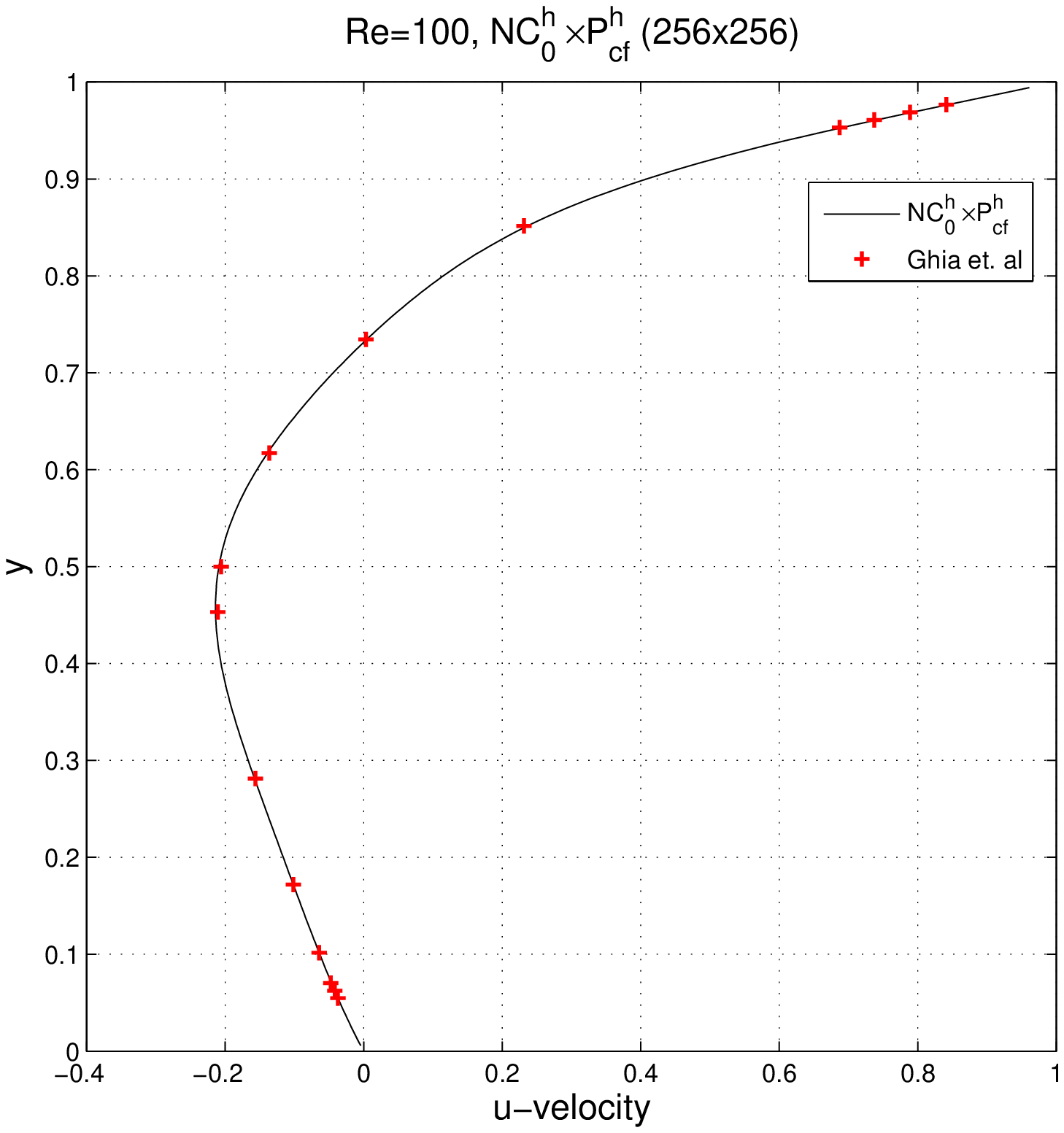}
\includegraphics[width=0.49\textwidth]{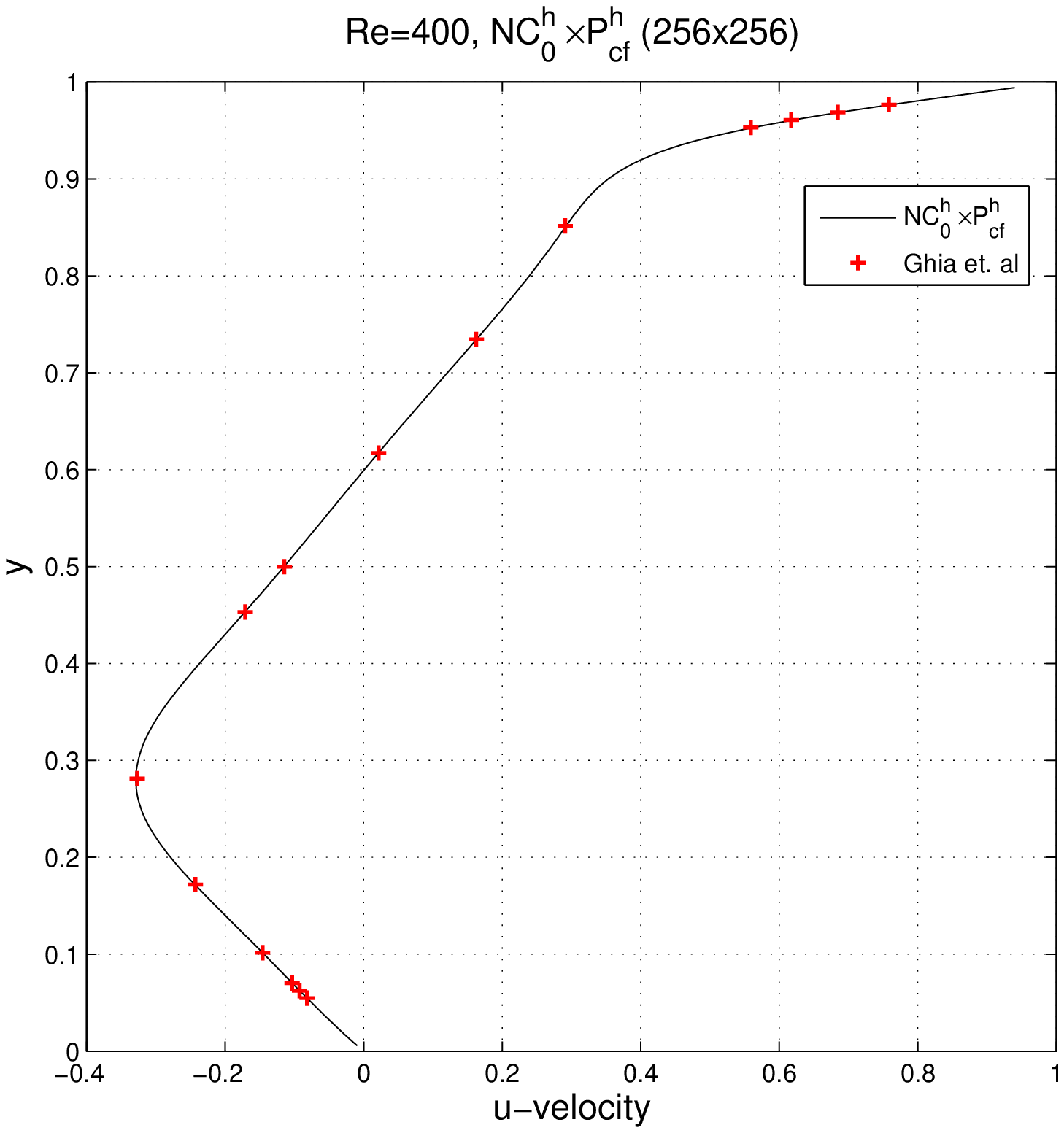}\\
\includegraphics[width=0.49\textwidth]{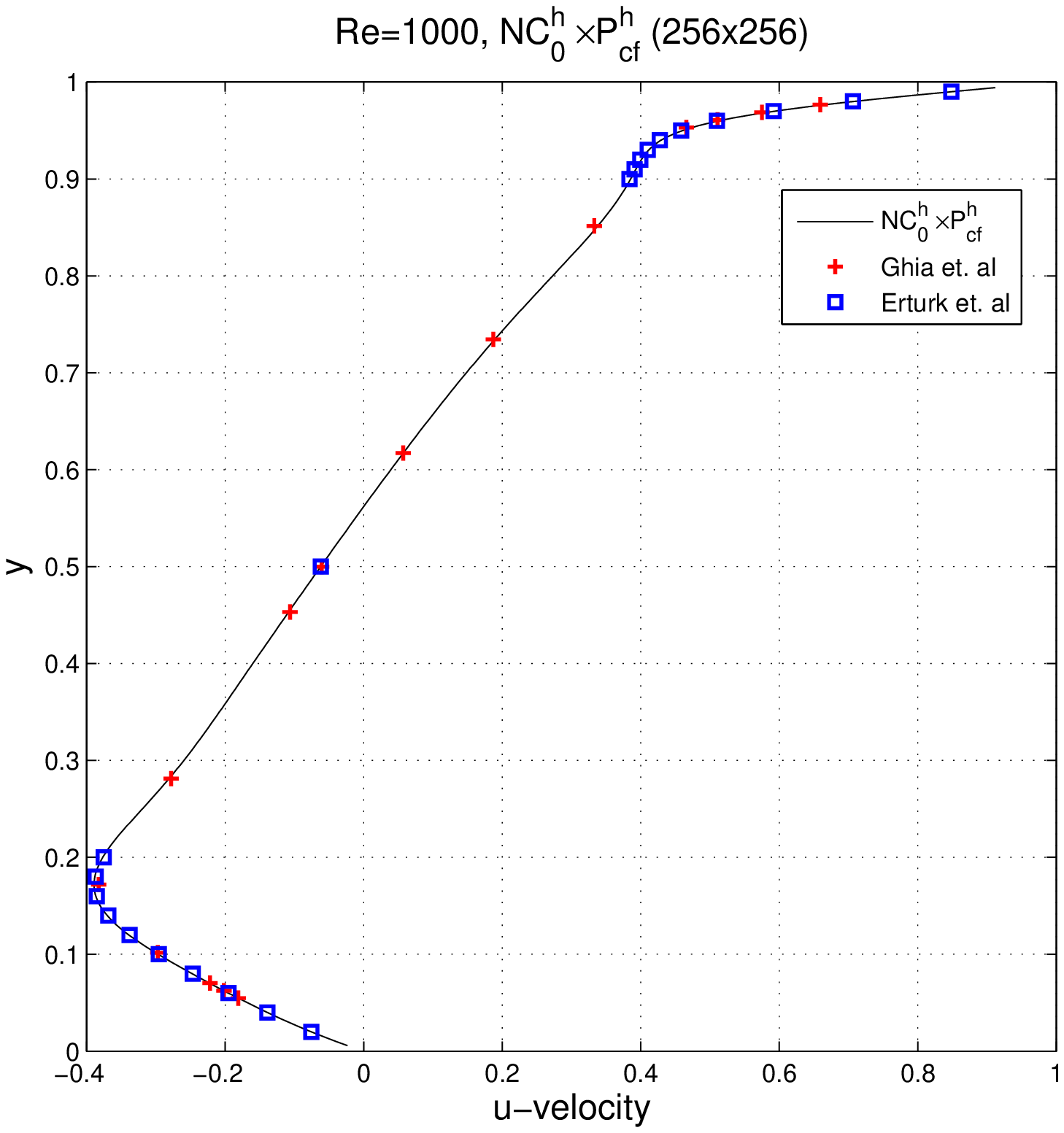}
\includegraphics[width=0.49\textwidth]{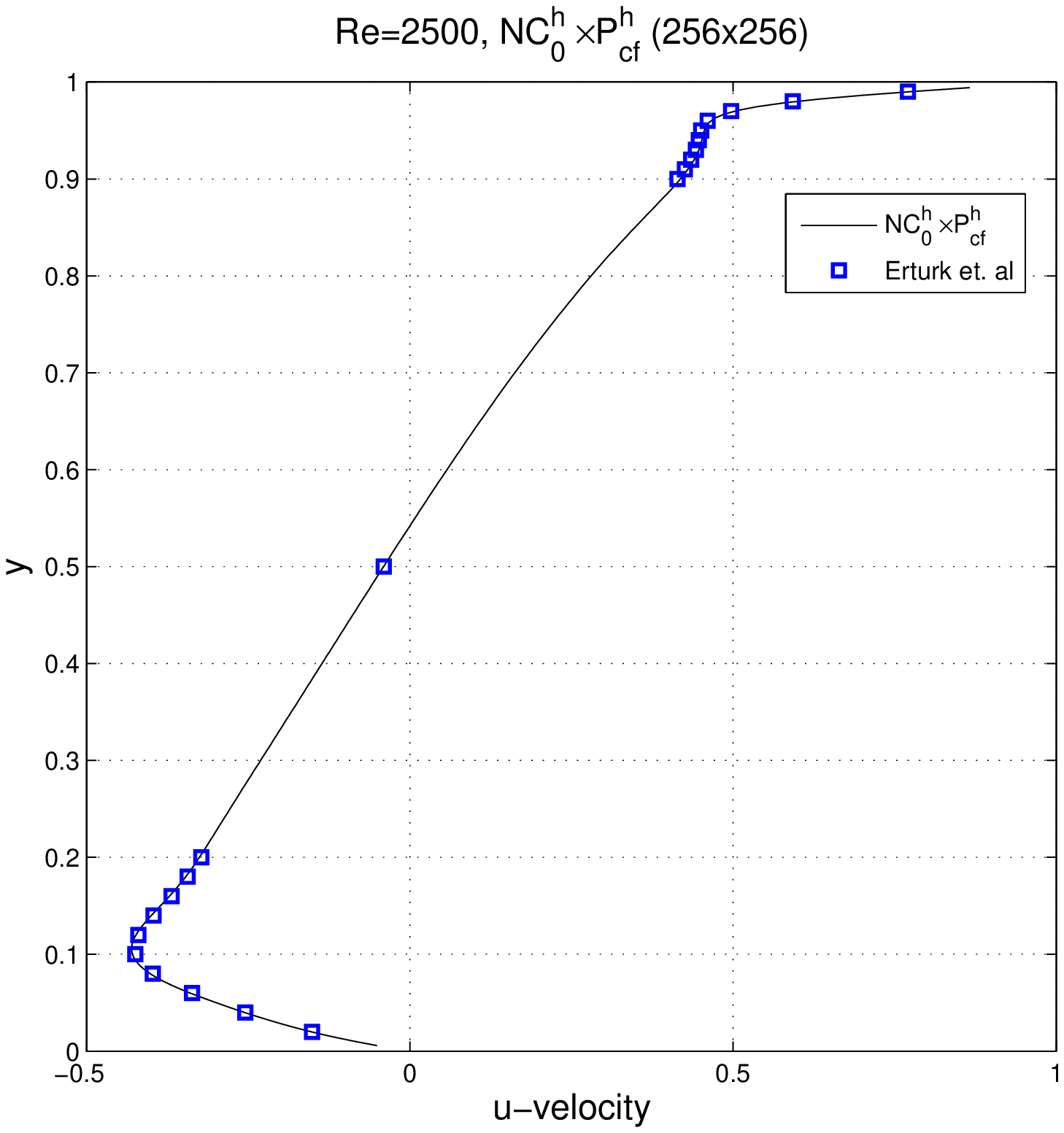}\\
\includegraphics[width=0.49\textwidth]{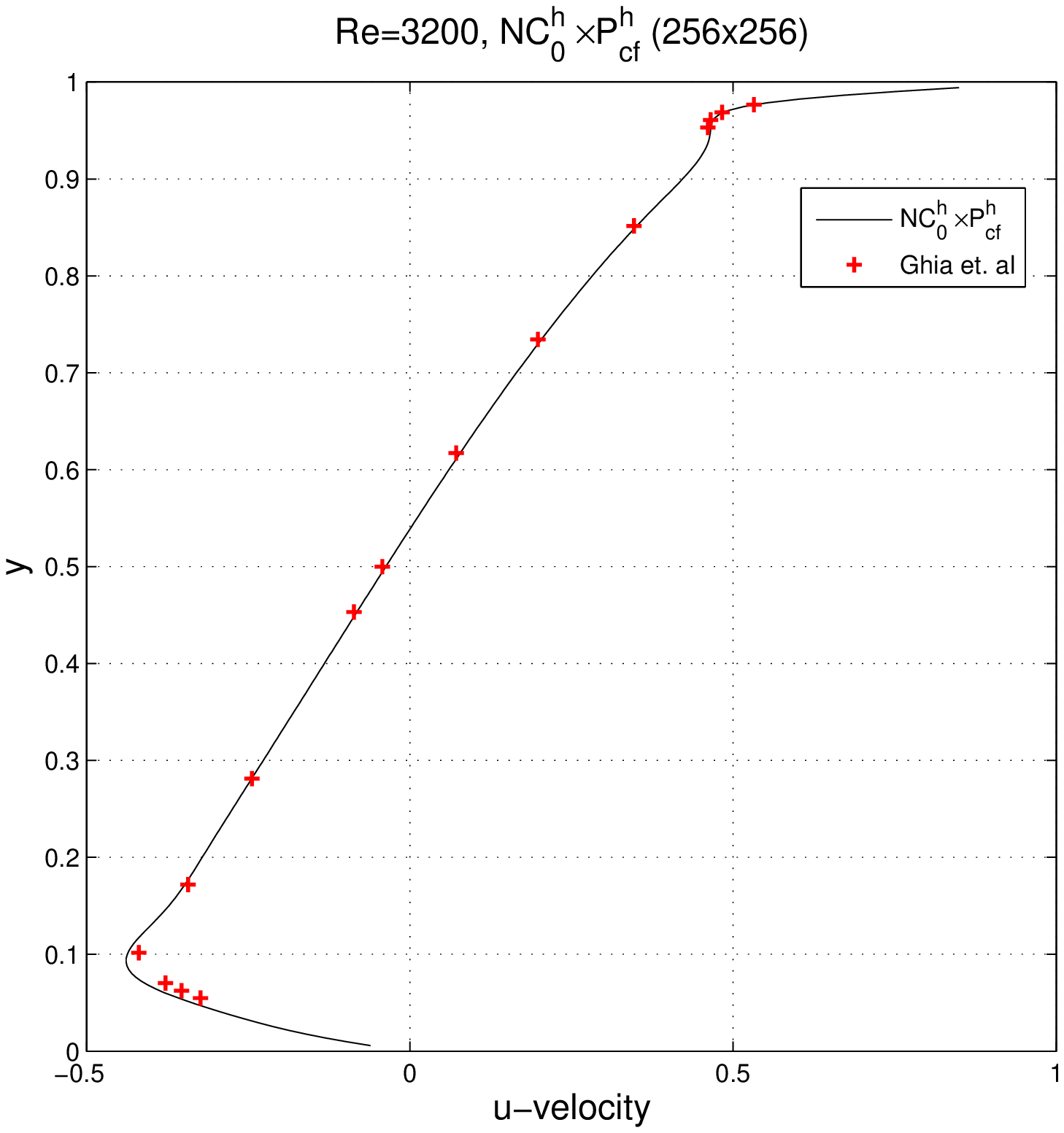}
\includegraphics[width=0.49\textwidth]{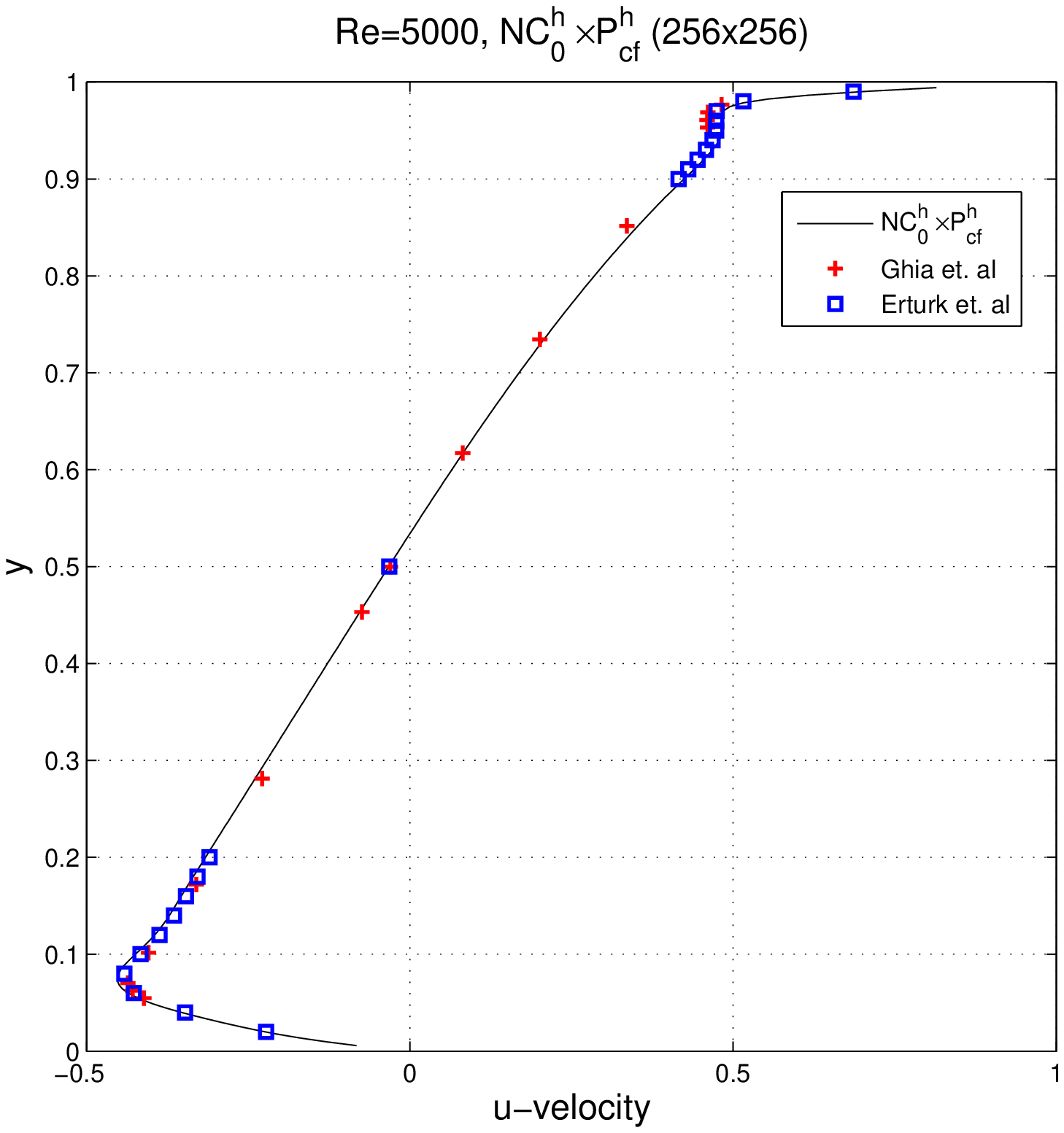}
\end{figure}

\begin{figure}
\caption{Profiles of $v$-velocity along the line $y=0.5$ computed by using the
$\vNChz \times \Pcf$ with unregularized boundary condition}
\label{vprofile}
\centering
\includegraphics[width=0.49\textwidth]{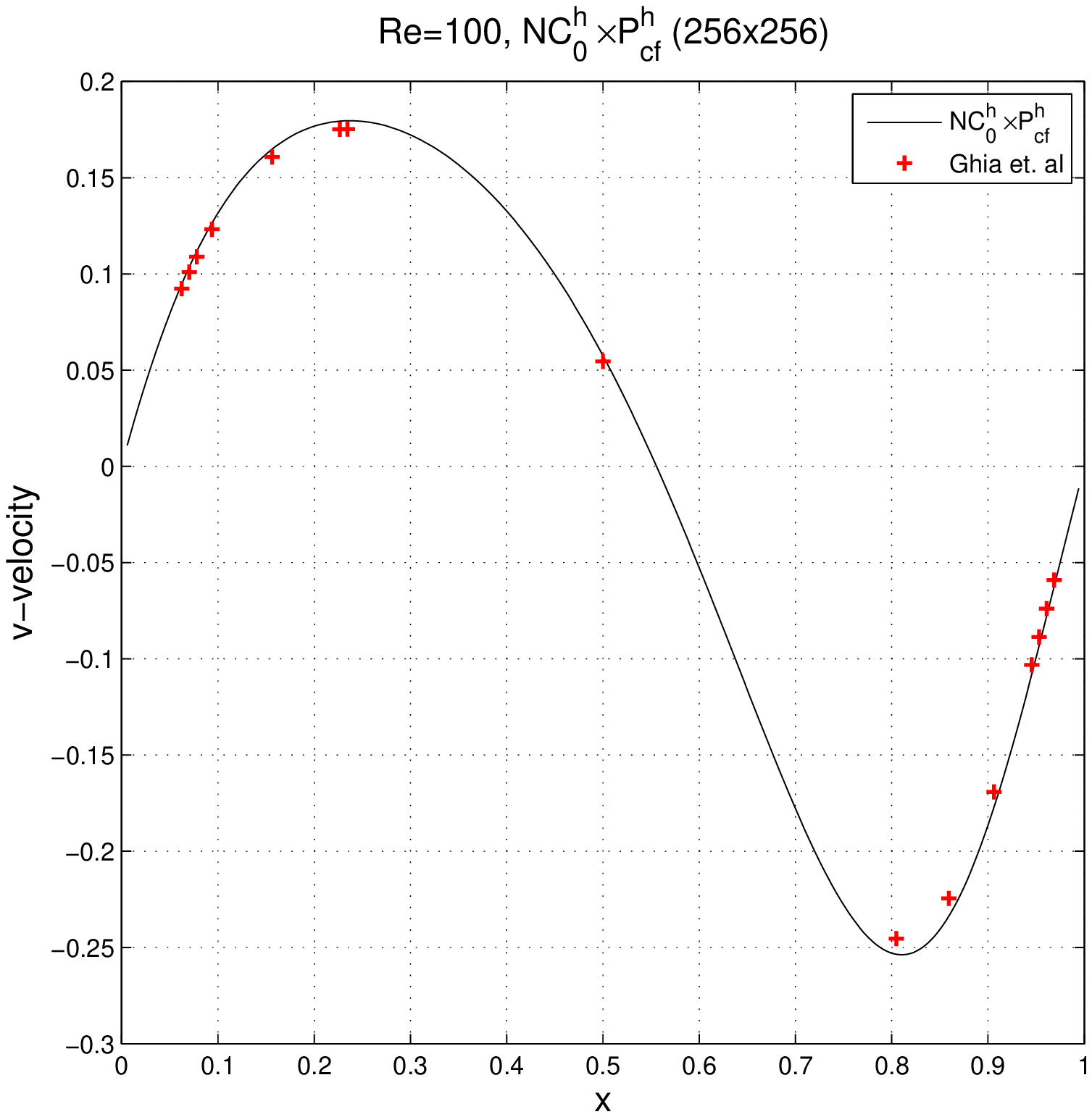}
\includegraphics[width=0.49\textwidth]{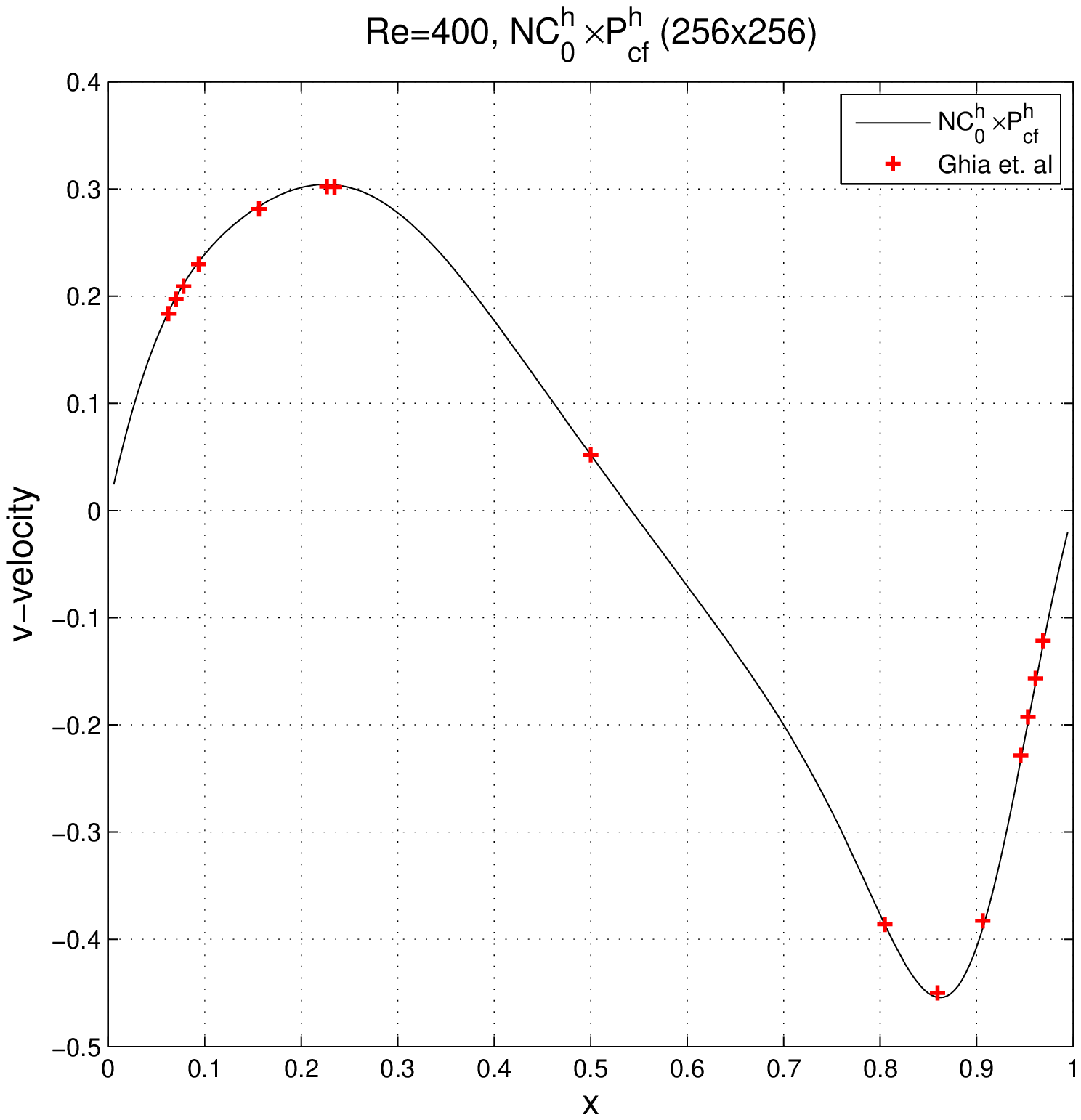}\\
\includegraphics[width=0.49\textwidth]{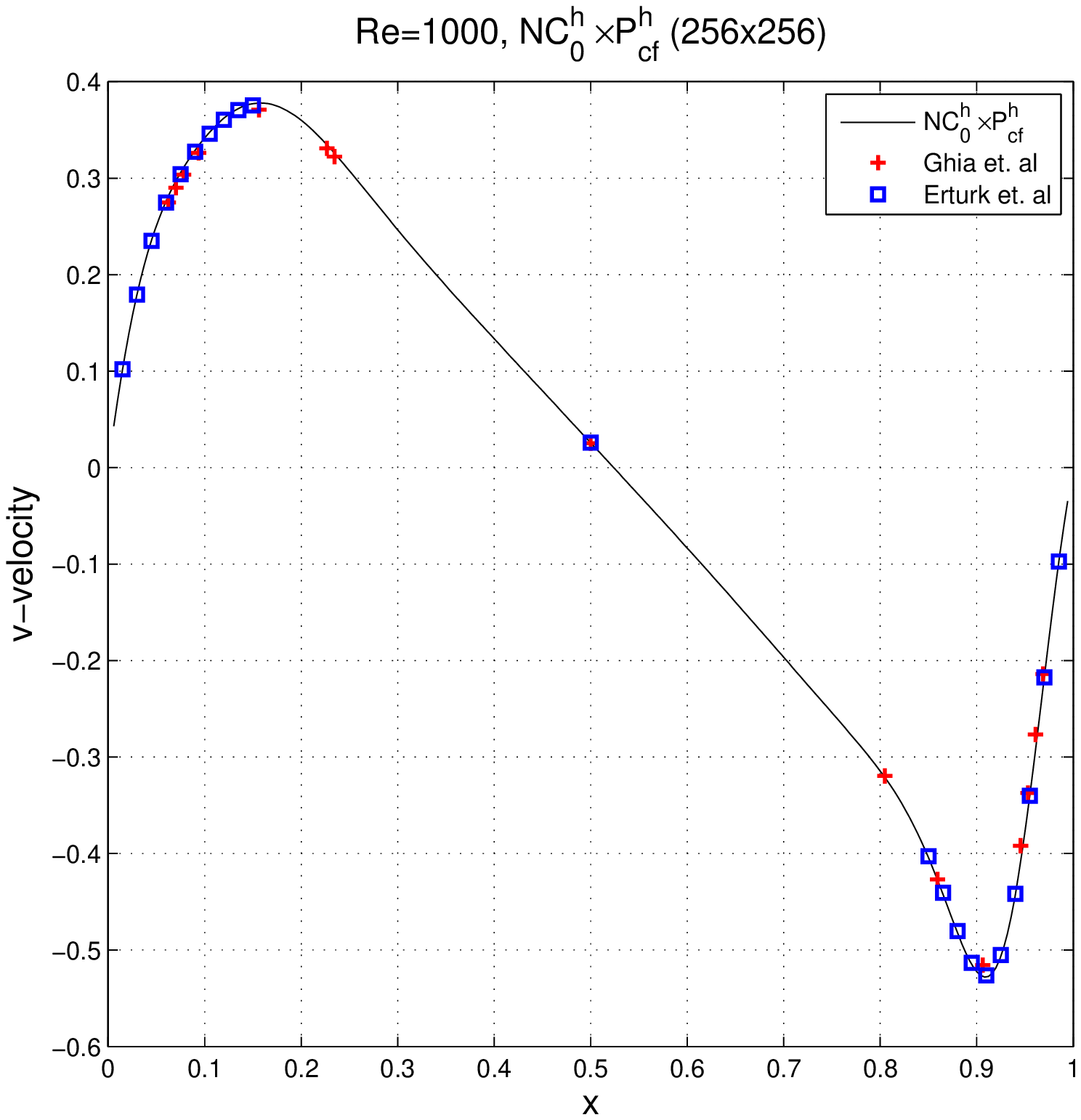}
\includegraphics[width=0.49\textwidth]{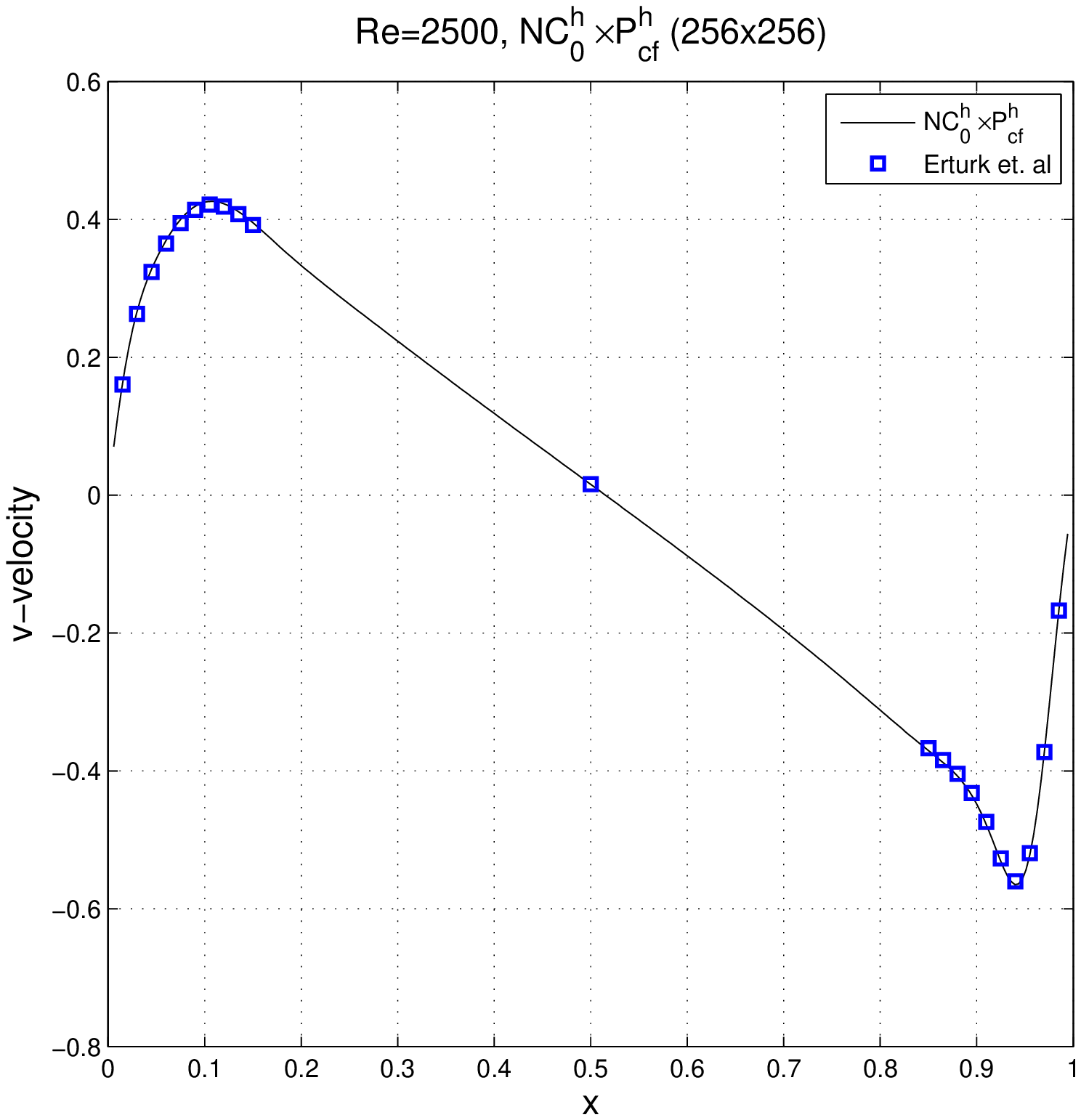}\\
\includegraphics[width=0.49\textwidth]{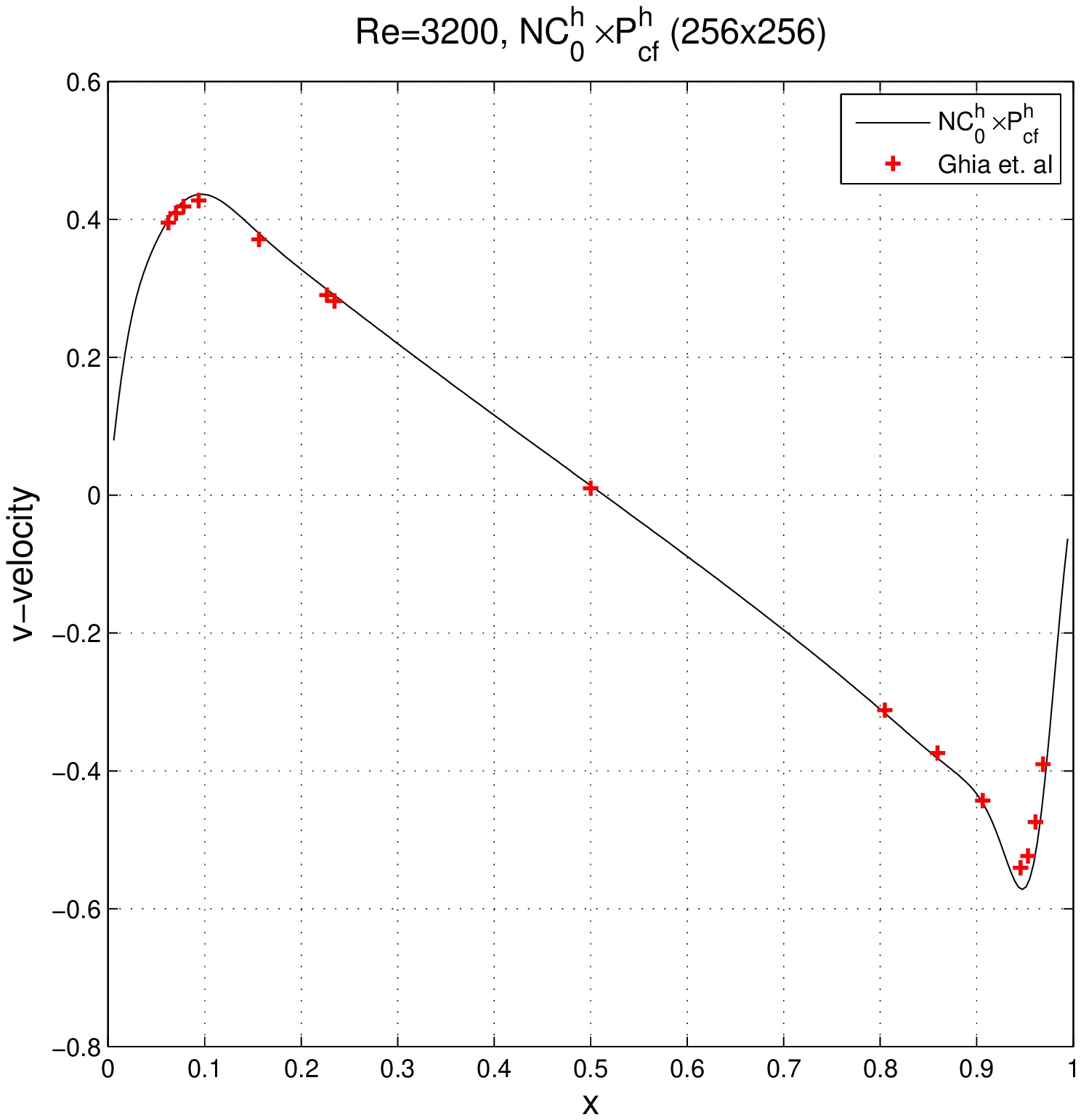}
\includegraphics[width=0.49\textwidth]{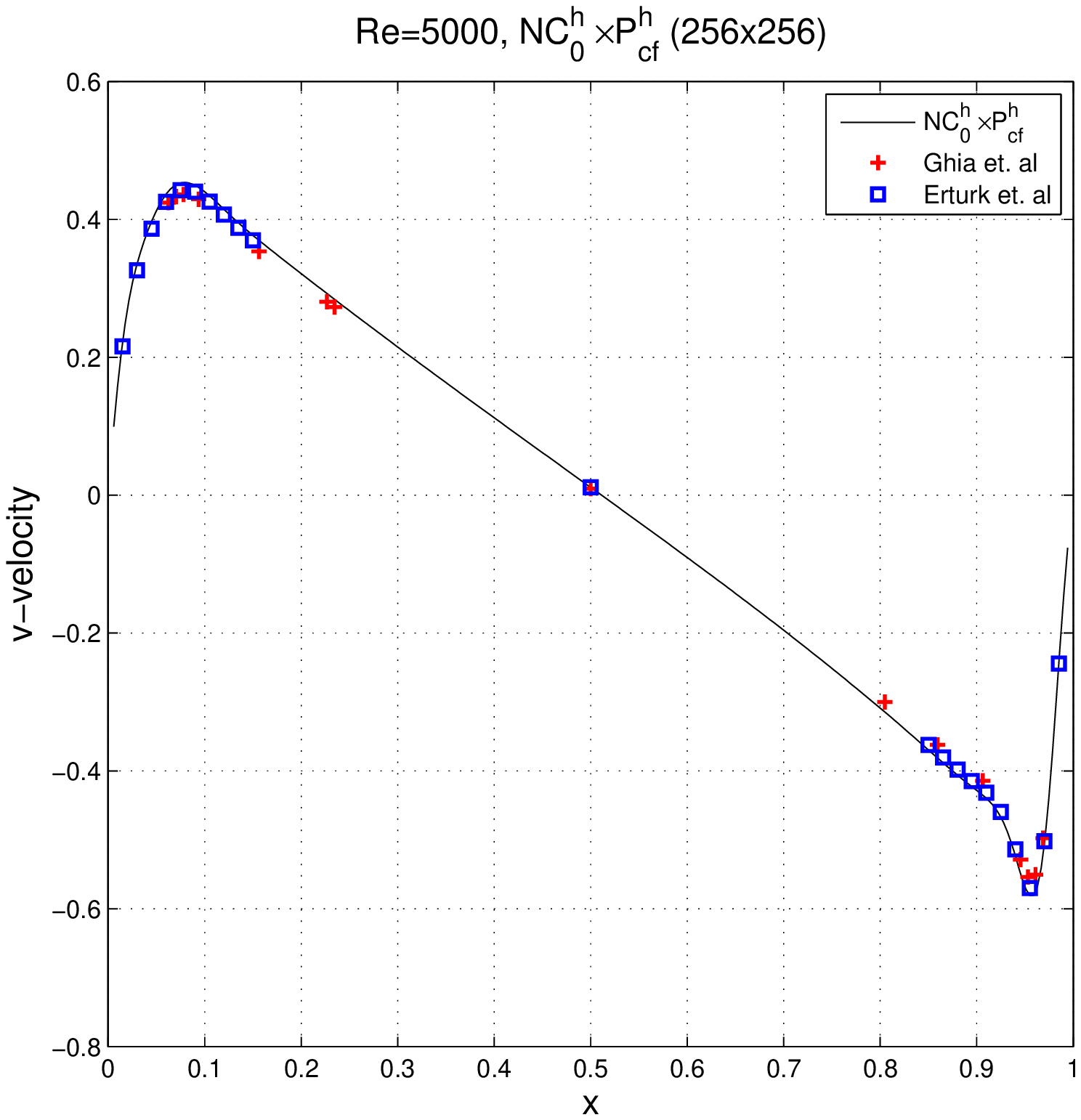}
\end{figure}

For the range of $Re=100, 400, 1000, 2500, 3200,$ and 5000,
Figs. \ref{uprofile} and \ref{vprofile} show the $u$-velocity profiles
along the line $x=0.5$ and the $v$-velocity profiles along the line $y=0.5$ computed
by using the $\vNChz \times \Pcf$ element with the boundary condition \eqref{cavityBD}
and comparison results with those by Erturk {\it et al.} \cite{erturk2005numerical} and
Ghia {\it et al.} \cite{ghia}. In each case, our velocity profiles
show a good agreement with their results.

Although the velocity profiles in Fig. \ref{fig:velocity-profile} and 
Figs. \ref{uprofile} and \ref{vprofile} seem to match quite well for $Re=1000$,
some of the actual numerical values differ in digits compared to those
reported in \cite{botella1998benchmark, bruneau20062d,
  guermond2012start}.
Hence, we compare the numerical values of the horizontal and vertical
components of the velocity in Tables 4 and 5
with the reference solutions from
\cite{botella1998benchmark, bruneau20062d, guermond2012start}. Our
numerical values, which were computed with $256\times 256$ meshes with the
lowest possible finite element $\vNChz \times \Pcf$, match with the
reference values mostly up to two digits, or with less than 1\% errors;
the numerical solutions, computed $512\times 512$ meshes
match with the
reference values mostly up to three digits, or with less than 0.1\% errors.

\begin{table}[htb!]\label{tab:horizontal-v}
\caption{Comparison of the horizontal components of the velocity
along the segment $y \in [0,1]$, $x=1/2$ at $Re=1000$}
\centering\begin{scriptsize}
\begin{tabular}{cccccc}
\hline
$y$ & \cite{botella1998benchmark} & \cite{bruneau20062d} &
\cite{guermond2012start} &
$\vNChz \times \Pcf (256 \times 256)$ &
$\vNChz \times \Pcf (512 \times 512)$ \\
\hline
 0.0000 &  0.0000000 &  0.00000 &  0.0000000 &  0.0000000 &  0.0000000\\
 0.0312 & -0.\emp{2279}225 &  NA      & -0.\emp{2279}177 & -0.\emp{22}74204 & -0.\emp{227}6650\\
 0.0391 & -0.\emp{2936}869 & -0.29330 & -0.\emp{2936}814 & -0.\emp{29}30076 & -0.\emp{293}3552\\
 0.0469 & -0.\emp{3553}213 &  NA      & -0.\emp{3553}154 & -0.\emp{35}45665 & -0.\emp{35}49485\\
 0.0547 & -0.\emp{4103}754 & -0.41018 & -0.\emp{4103}691 & -0.\emp{40}96654 & -0.\emp{410}0002\\
 0.0937 & -0.\emp{5264}392 &  NA      & -0.\emp{5264}320 & -0.\emp{52}71749 & -0.\emp{526}4518\\
 0.1406 & -0.\emp{4264}545 & -0.42645 & -0.\emp{4264}492 & -0.\emp{42}76315 & -0.\emp{426}5356\\
 0.1953 & -0.\emp{3202}137 &  NA      & -0.\emp{3202}068 & -0.\emp{32}09943 & -0.\emp{320}0577\\
 0.5000 &  0.\emp{0257}995 &  0.02580 &  0.\emp{0257}987 &  0.\emp{02}56839 &  0.\emp{025}7175\\
 0.7656 &  0.\emp{3253}592 &  NA      &  0.\emp{3253}529 &  0.\emp{32}59697 &  0.\emp{325}2217\\
 0.7734 &  0.\emp{3339}924 &  0.33398 &  0.\emp{3339}860 &  0.\emp{33}46373 &  0.\emp{333}8694\\
 0.8437 &  0.\emp{3769}189 &  NA      &  0.\emp{3769}119 &  0.\emp{37}78450 &  0.\emp{376}9140\\
 0.9062 &  0.\emp{3330}442 &  0.33290 &  0.\emp{3330}381 &  0.\emp{33}39829 &  0.\emp{333}1021\\
 0.9219 &  0.\emp{3099}097 &  NA      &  0.\emp{3099}041 &  0.\emp{31}08006 &  0.\emp{309}9725\\
 0.9297 &  0.\emp{2962}703 &  0.29622 &  0.\emp{2962}650 &  0.\emp{2}971221 &  0.\emp{296}3312\\
 0.9375 &  0.\emp{2807}056 &  NA      &  0.\emp{2807}005 &  0.\emp{28}15029 &  0.\emp{280}7605\\
 1.0000 &  0.0000000 &  0.00000 &  0.0000000 &  0.0000000 &  0.0000000\\
\hline
\end{tabular}\end{scriptsize}
\end{table}

\begin{table}[htb!]\label{tab:vertical-u}
\caption{Comparison of the vertical components of the velocity
along the segment $x \in [0,1]$, $y=1/2$ at $Re=1000$}

\centering\begin{scriptsize}
\begin{tabular}{cccccc}
\hline
$x$ & \cite{botella1998benchmark} & \cite{bruneau20062d} &
\cite{guermond2012start} &
$\vNChz \times \Pcf (256 \times 256)$ &
$\vNChz \times \Pcf (512 \times 512)$ \\
\hline
 1.0000 & -1.0000000 & -1.00000 & -1.0000000 & -1.0000000 & -1.0000000\\
 0.9766 & -0.\emp{6644}227 & NA       & -0.\emp{6644}194 & -0.\emp{66}66343 & -0.\emp{664}8562\\
 0.9688 & -0.\emp{5808}359 & -0.58031 & -0.\emp{5808}318 & -0.\emp{58}31751 & -0.\emp{58}12660\\
 0.9609 & -0.\emp{5169}277 & NA       & -0.\emp{5169}214 & -0.\emp{51}90905 & -0.\emp{51}72781\\
 0.9531 & -0.\emp{4723}329 & -0.47239 & -0.\emp{4723}260 & -0.\emp{47}41970 & -0.\emp{472}5743\\
 0.8516 & -0.\emp{3372}212 & NA       & -0.\emp{3372}128 & -0.\emp{33}80993 & -0.\emp{337}0508\\
 0.7344 & -0.\emp{1886}747 & -0.18861 & -0.\emp{1886}680 & -0.\emp{18}90994 & -0.\emp{188}4232\\
 0.6172 & -0.\emp{0570}178 & NA       & -0.\emp{0570}151 & -0.\emp{05}70951 & -0.\emp{056}9011\\
 0.5000 &  0.\emp{0620}561 &  0.06205 &  0.\emp{0620}535 &  0.\emp{06}22962 & -0.\emp{061}9466\\
 0.4531 &  0.\emp{1081}999 & NA       &  0.\emp{1081}955 &  0.\emp{10}85611 &  0.\emp{108}0176\\
 0.2813 &  0.\emp{2803}696 &  0.28040 &  0.\emp{2803}632 &  0.\emp{28}11184 &  0.\emp{280}2013\\
 0.1719 &  0.\emp{3885}691 & NA       &  0.\emp{3885}624 &  0.\emp{38}94565 &  0.\emp{388}5914\\
 0.1016 &  0.\emp{3004}561 &  0.30029 &  0.\emp{3004}504 &  0.\emp{30}06758 &  0.\emp{300}4357\\
 0.0703 &  0.\emp{2228}955 & NA       &  0.\emp{2228}928 &  0.\emp{22}28075 &  0.\emp{222}8534\\
 0.0625 &  0.\emp{2023}300 &  0.20227 &  0.\emp{2023}277 &  0.\emp{20}21815 &  0.\emp{202}2834\\
 0.0547 &  0.\emp{1812}881 & NA       &  0.\emp{1812}863 &  0.\emp{18}10885 &  0.\emp{181}2376\\
 0.0000 &  0.0000000 &  0.00000 &  0.0000000 &  0.0000000 &  0.0000000\\
\hline
\end{tabular}\end{scriptsize}
\end{table}

\begin{figure}
\caption{Streamline computed by using the $\vNChz \times \Pcf$ element.}
\label{fig:streamline}
\centering
\includegraphics[width=0.49\textwidth]{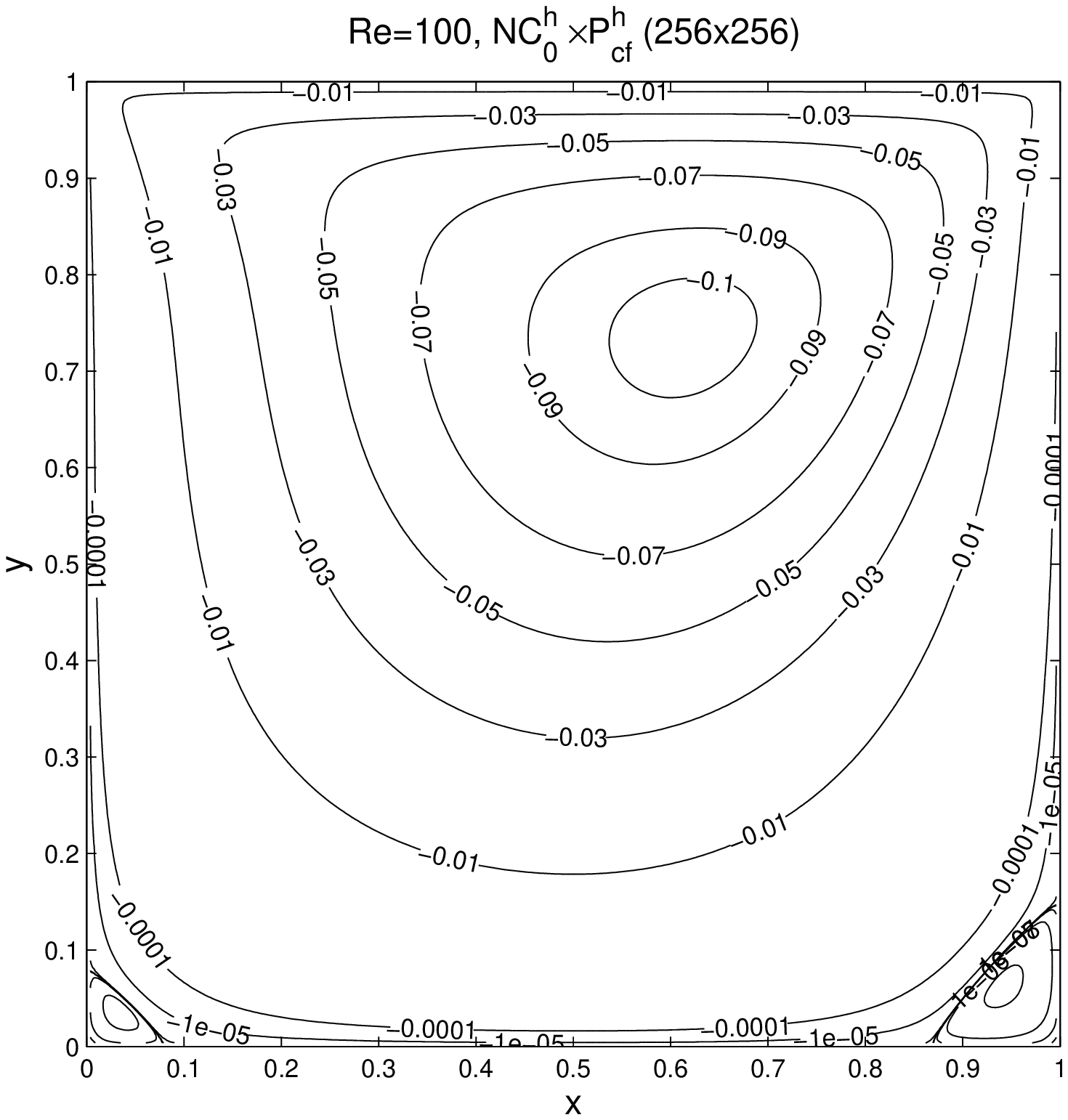}
\includegraphics[width=0.49\textwidth]{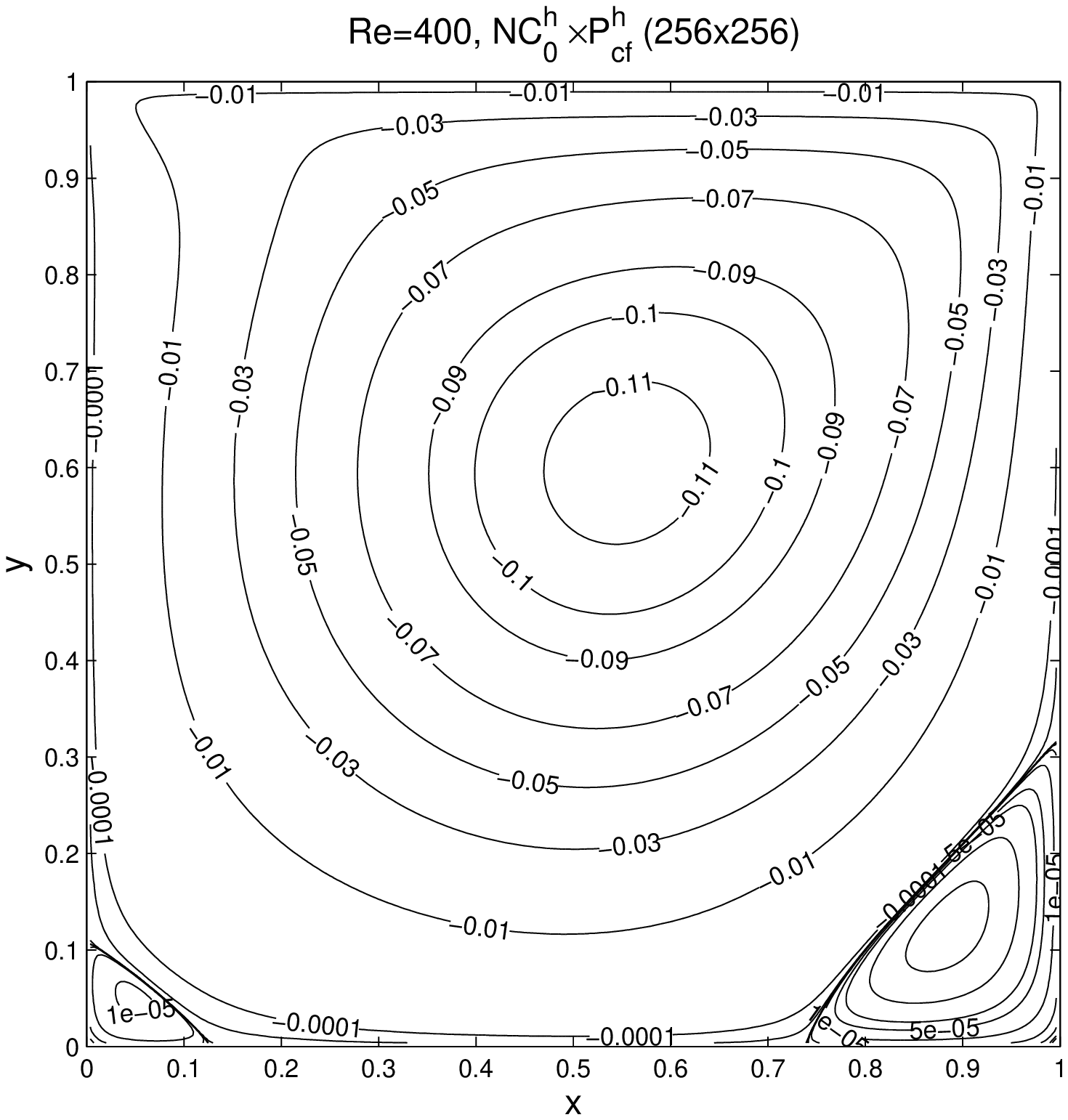}\\
\includegraphics[width=0.49\textwidth]{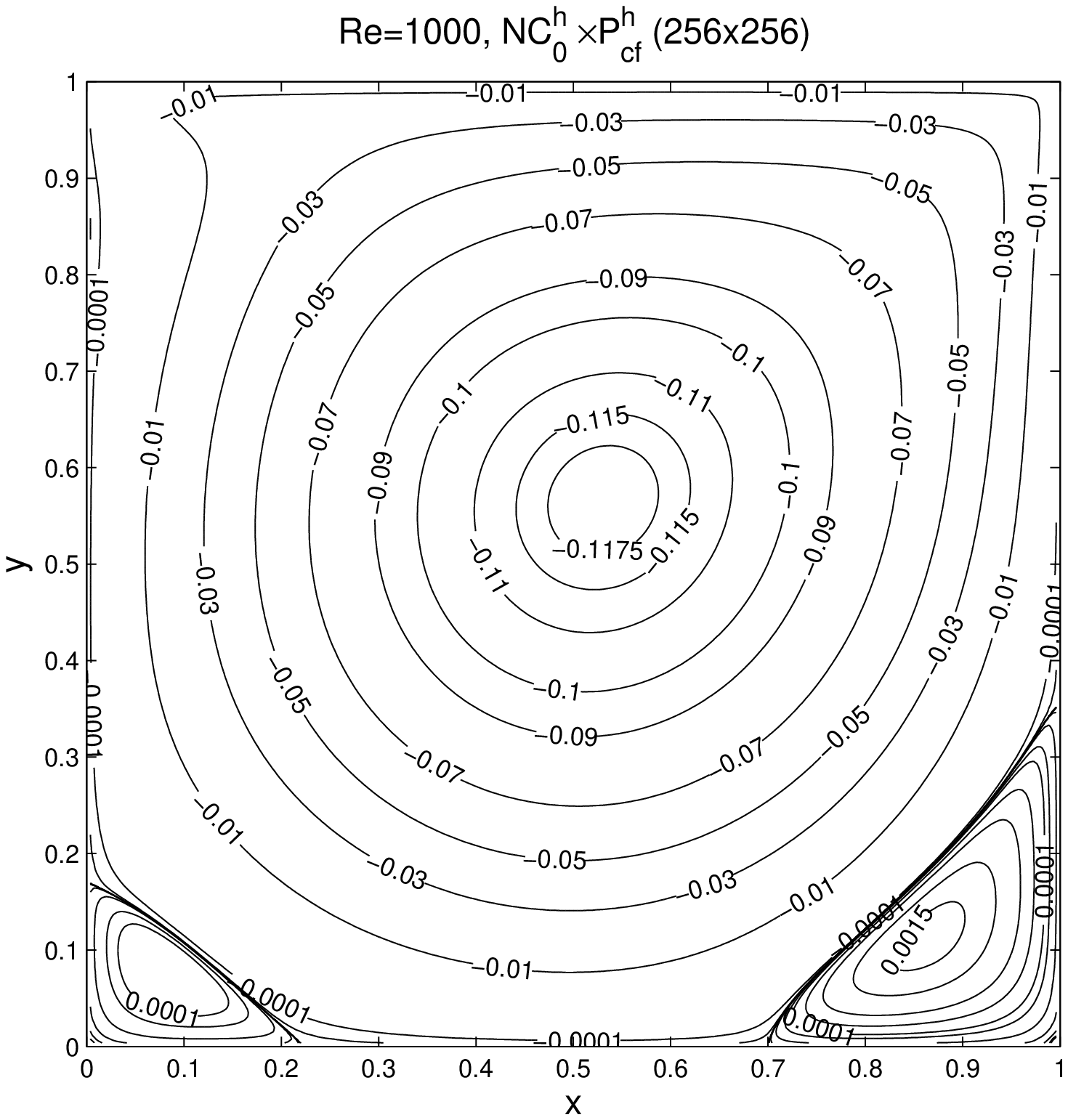}
\includegraphics[width=0.49\textwidth]{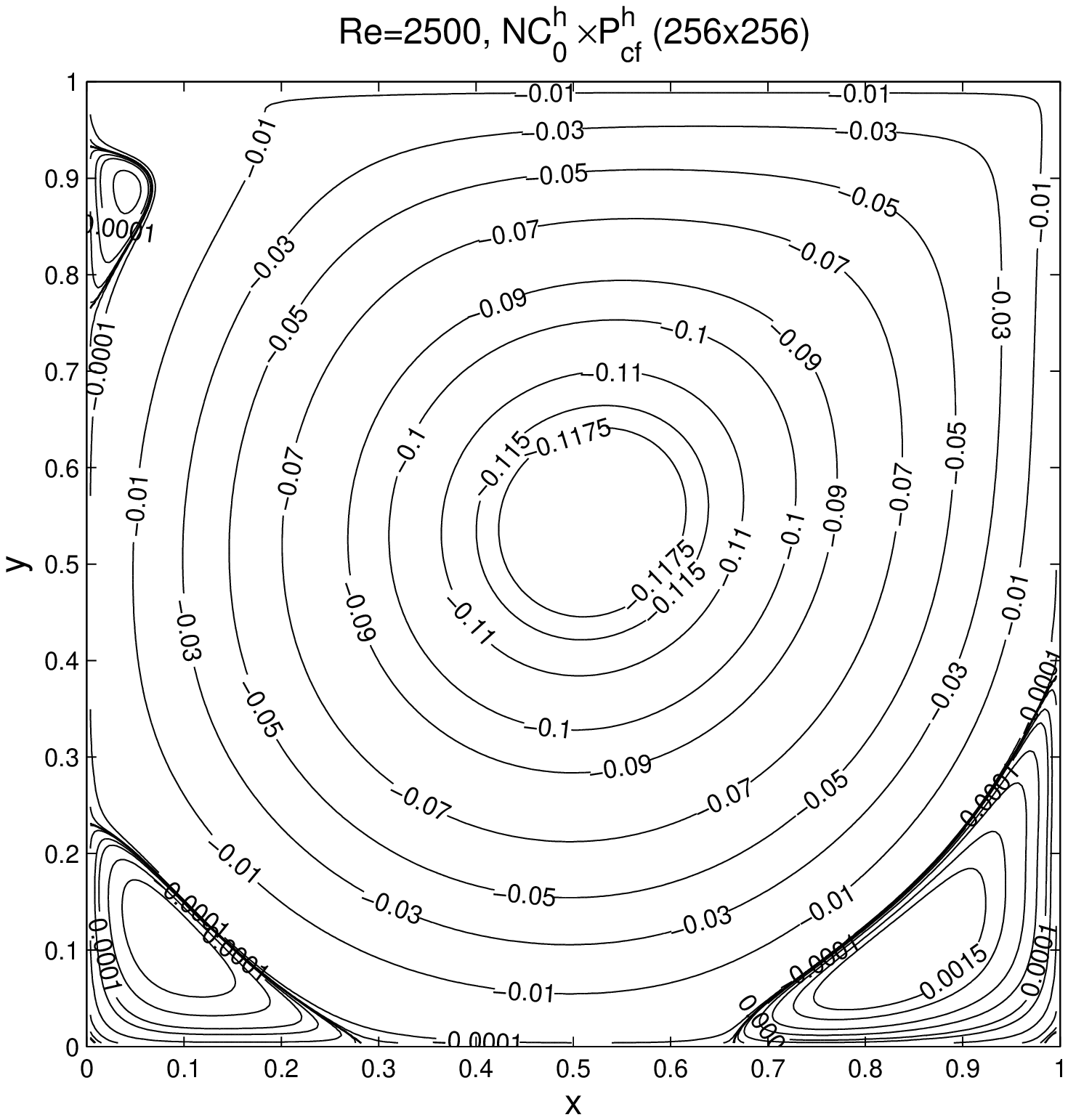}\\
\includegraphics[width=0.49\textwidth]{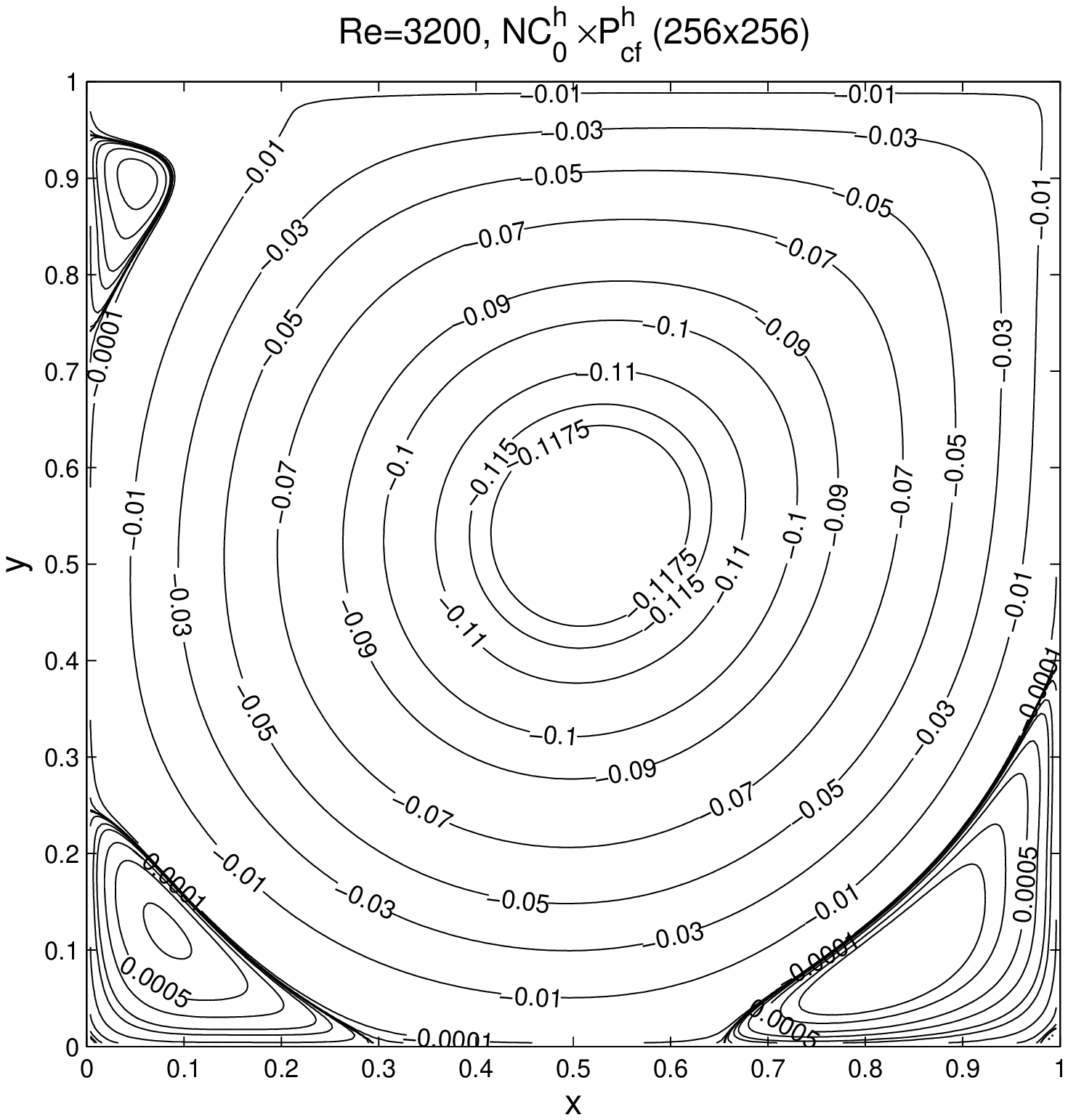}
\includegraphics[width=0.49\textwidth]{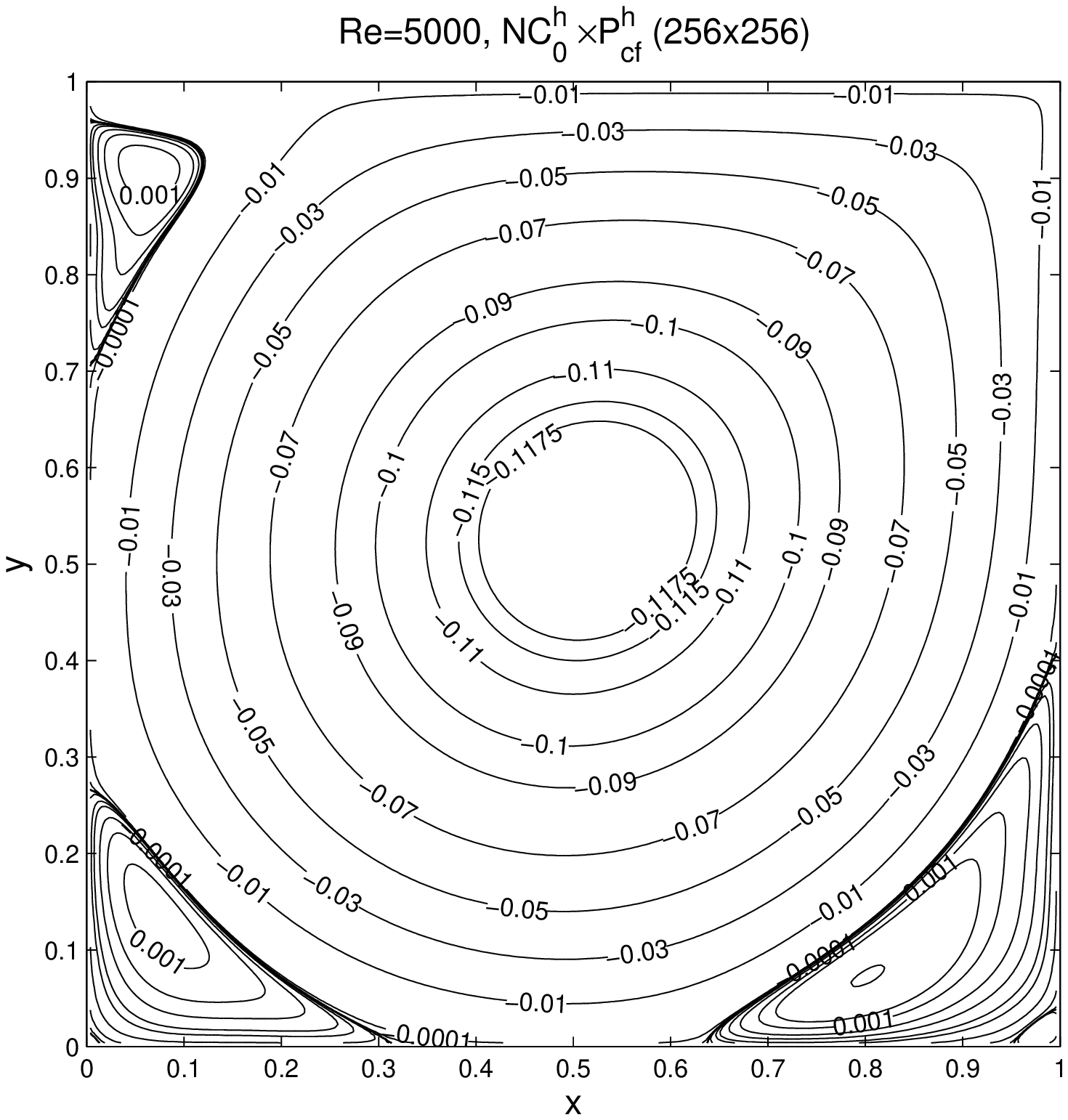}
\end{figure}

The computed streamlines are presented in \figref{fig:streamline}.
One can observe count-rotating secondary vortices at the bottom left
and right corners of the square cavity. Bottom left and right vortices grow
in size as Reynolds number increases and the secondary vortex at the top
left corner of the square cavity develops as Reynolds number increases.

\begin{figure}
\caption{Contours of vorticity computed by using the $\vNChz \times \Pcf$
element.}
\label{fig:vorticity}
\centering
\includegraphics[width=0.49\textwidth]{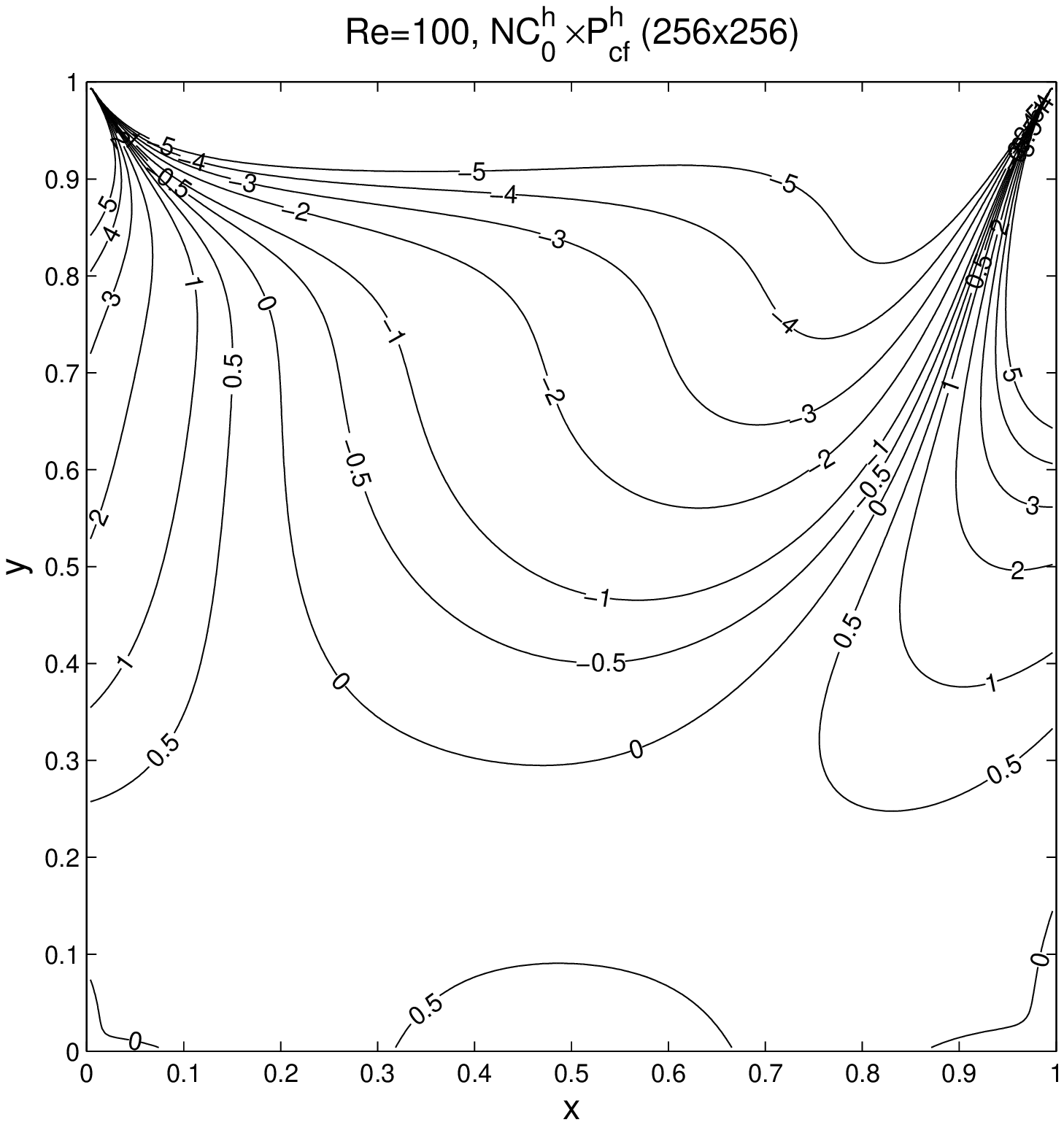}
\includegraphics[width=0.49\textwidth]{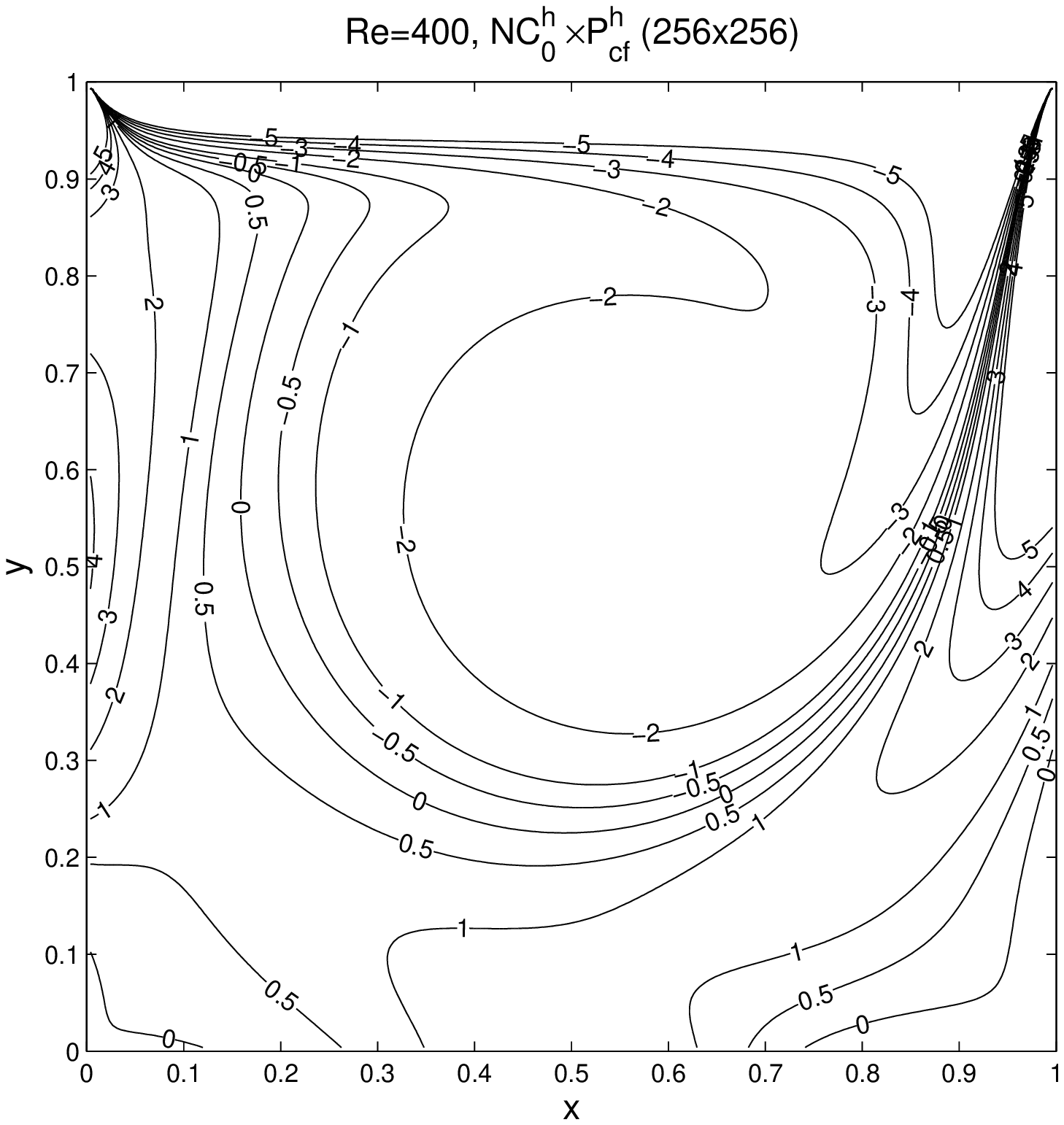}\\
\includegraphics[width=0.49\textwidth]{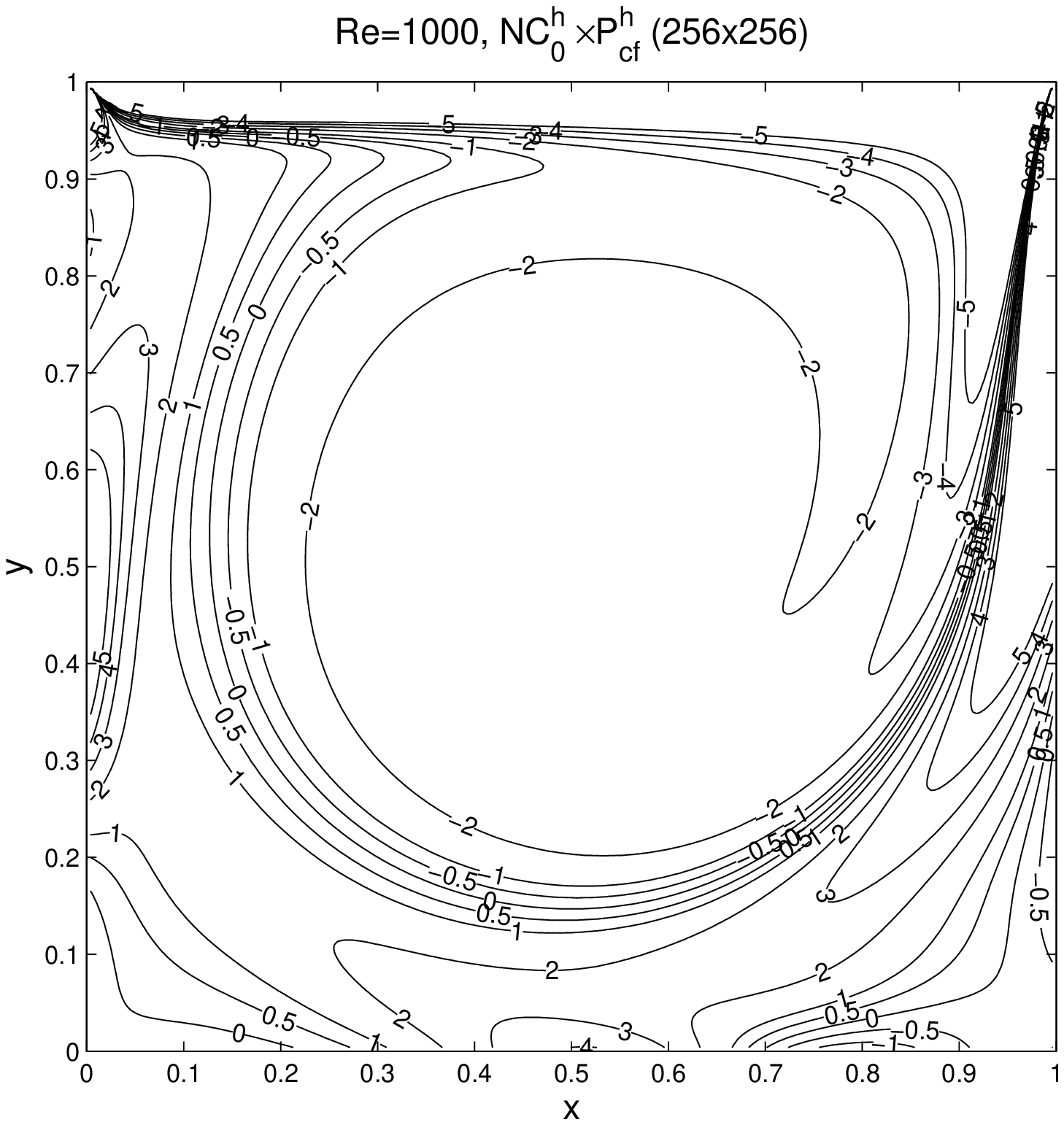}
\includegraphics[width=0.49\textwidth]{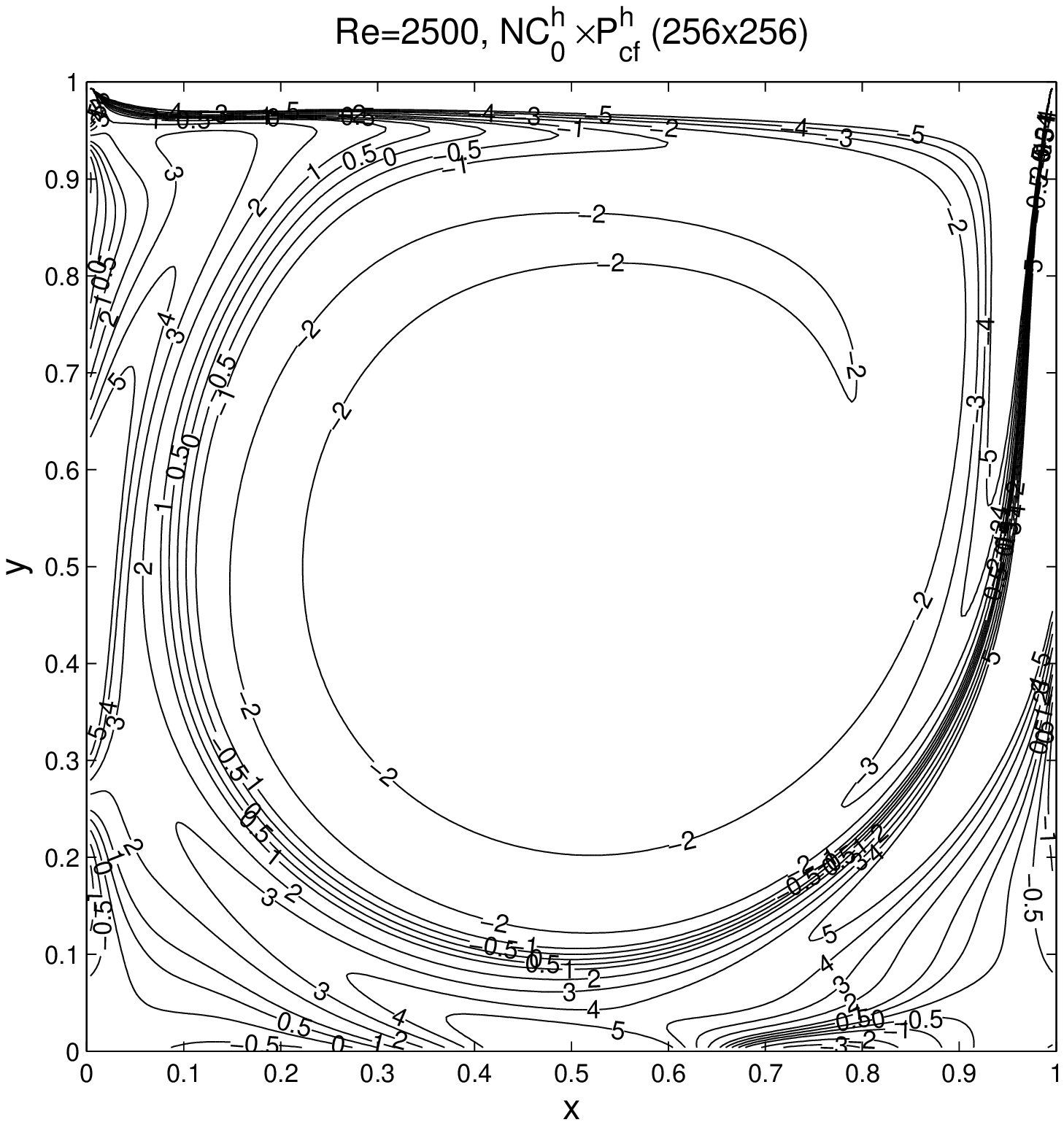}\\
\includegraphics[width=0.49\textwidth]{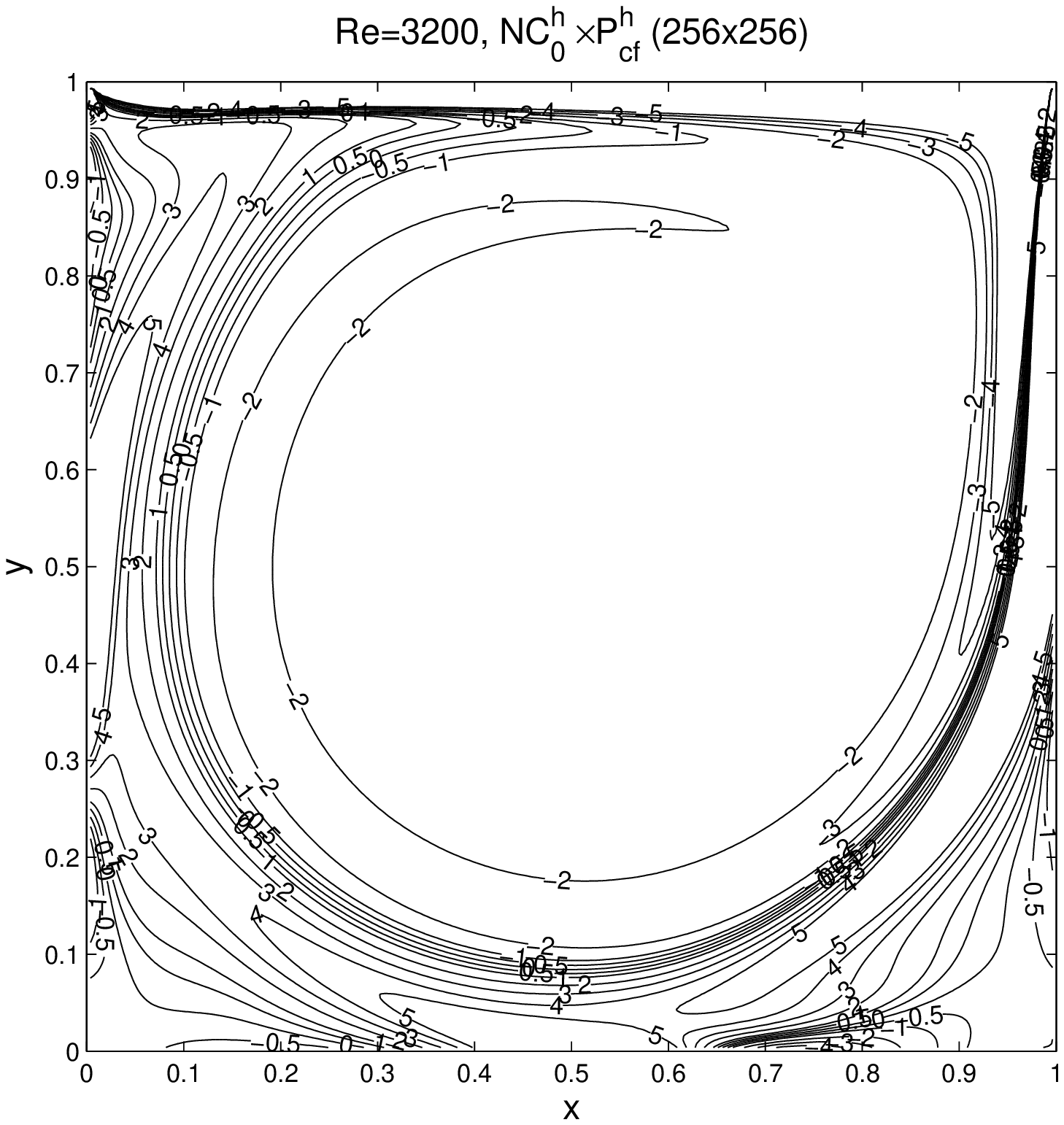}
\includegraphics[width=0.49\textwidth]{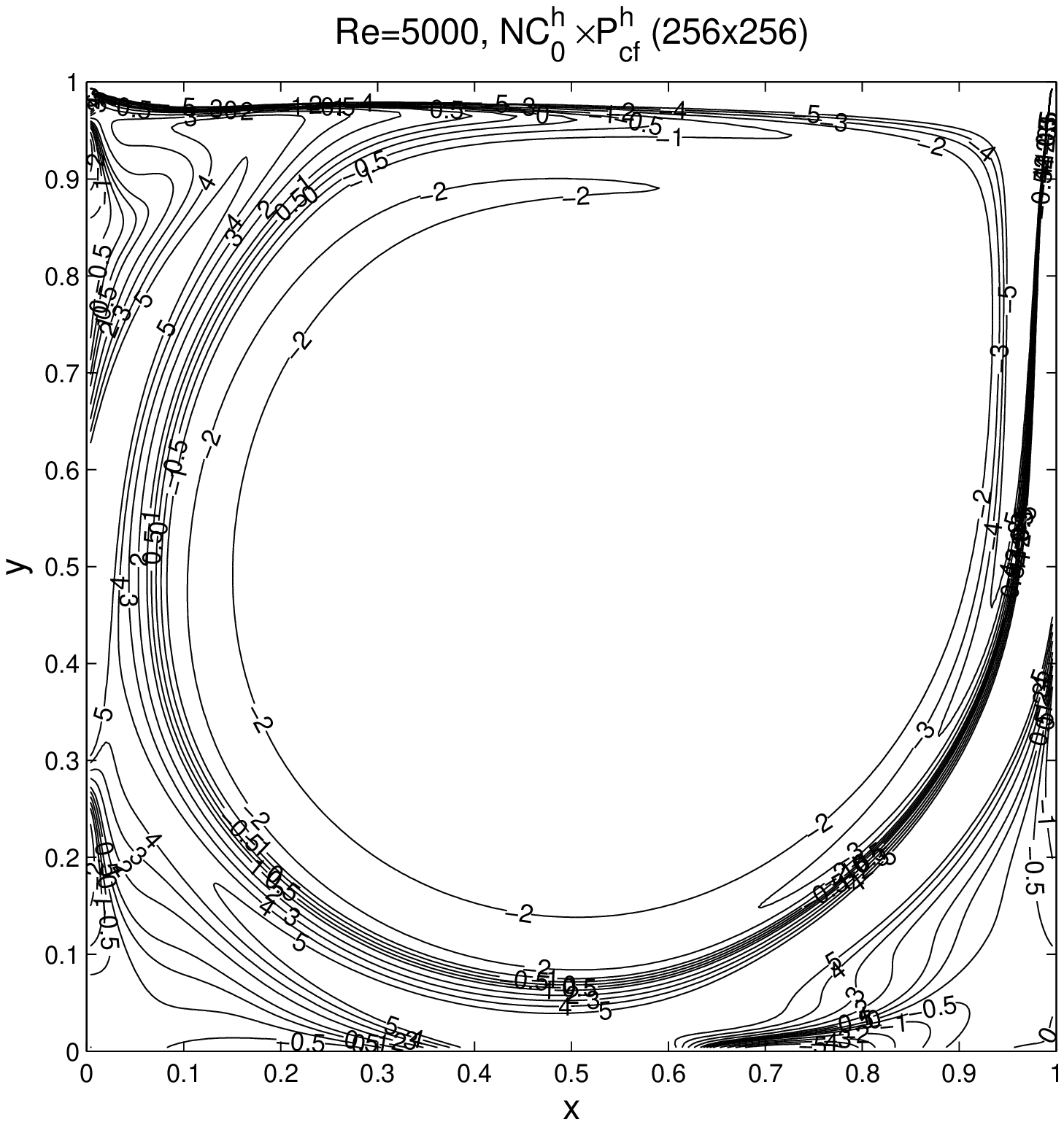}
\end{figure}

The vorticity contours are presented in \figref{fig:vorticity}.
We observe that the gradient in vorticity is negligible
in the center of cavity and the region of very low gradient
in vorticity grows as Reynolds number increases.

In \tabref{propvortex}, we present the location of the center of the primary
vortex, the stream function $\psi$, and vorticity $\omega$ at vortex center.
These data are calculated for $100 \leq \text{Re} \leq 5000$; for comparison,
available data from the literatures are also given.
The values of the stream function $\psi$ and vorticity $\omega$ are recorded
at the center of meshes. 
The locations of primary vortices computed by using the $\vNChz \times \Pcf$ element
differ from the other results by about $0.002$ which is half the mesh size
$1/h \approx 0.0039$, Our numerical solutions computed by using both 
$\vNChz \times \Pcf$ and $\vcQqhz \times \dcQlhz$ elements exhibit a good agreement
with the literature data except in the case of the $\vcQqhz \times \dcQlhz$
element with leaky cavity boundary condition \eqref{leaky} applied.
$\psi$ and $\omega$ for the $\vNChz \times \Pcf$ element with \eqref{cavityBD}
are similar to those in the literature
\cite{erturk2005numerical, ghia, hna09, sahin2003}.
In addition, \tabref{secondvortex} summarizes data on the strengths and the locations of
secondary vortices in the bottom left and right corners, and in the top left
corner. We observe that secondary vortices appear stronger as the Reynolds number increases.
 
\begin{table}
\caption{Computed primary vortex;
the values of stream function ($\psi$), vorticity ($\omega$),
and location $(x,y)$.} \label{propvortex}
\centering
\begin{scriptsize}
\begin{tabular}{c|cccccc}
\hline
Re & FEM & Grid & $\psi_{\min}$ & $\omega$ & $(x,y)$ & BC \\ \hline
 & $\vNChz \times \Pcf$ & $256\times 256$ &
-0.103531 & 3.16206 & (0.6152,0.7363) & \eqref{cavityBD} \\
& $\vcQqhz \times \dcQlhz$ & $128\times 128$ &
-0.103519 & 3.18101 & (0.6172,0.7383) & \eqref{watertight} \\
& $\vcQqhz \times \dcQlhz$ & $128\times 128$ &
-0.102872 & 3.15485 & (0.6172,0.7383) & \eqref{leaky} \\
\rb{100} & \cite{ghia} & $129\times 129$ &
-0.103423 & 3.16646 & (0.6172,0.7344) & - \\
& \cite{hna09} & $128\times 128$ &
-0.103435 & - & (0.6172,0.7344) & \cite{hna09} \\
& \cite{sahin2003} & $257\times 257$ &
-0.103471 & 3.1655 & (0.6189,0.7400) & \cite{sahin2003} \\
\hline 
 & $\vNChz \times \Pcf$ & $256\times 256$ &
-0.114071 & 2.29821 & (0.5527,0.6035) & \eqref{cavityBD} \\
& $\vcQqhz \times \dcQlhz$ & $128\times 128$ &
-0.113990 & 2.29476 & (0.5547,0.6055) & \eqref{watertight} \\
& $\vcQqhz \times \dcQlhz$ & $128\times 128$ &
-0.111900 & 2.26041 & (0.5547,0.6055) & \eqref{leaky} \\
\rb{400} & \cite{ghia} & $257\times 257$ &
-0.113909 & 2.29469 & (0.5547,0.6055) & - \\
& \cite{hna09} & $128\times 128$ &
-0.113909 & - & (0.5547,0.6094) & \cite{hna09} \\
& \cite{sahin2003} & $257\times 257$ &
-0.113897 & 2.2950 & (0.5536,0.6075) & \cite{sahin2003} \\
\hline
& $\vNChz \times \Pcf$ & $256\times 256$ &
-0.119186 & 2.07216 & (0.5293,0.5645) & \eqref{cavityBD} \\
& $\vcQqhz \times \dcQlhz$ & $128\times 128$ &
-0.118941 & 2.06779 & (0.5313,0.5664) & \eqref{watertight} \\
& $\vcQqhz \times \dcQlhz$ & $128\times 128$ &
-0.115376 & 2.00941 & (0.5313,0.5664) & \eqref{leaky} \\
1000 & \cite{erturk2005numerical} & $601 \times 601$ &
-0.118781 & 2.06553 & (0.5300,0.5650) & - \\
& \cite{ghia} & $257\times 257$ &
-0.117929 & 2.04968 & (0.5313,0.5625) & - \\
& \cite{hna09} & $128\times 128$ &
-0.119173 & - & (0.5313,0.5625) & \cite{hna09} \\
& \cite{sahin2003} & $257\times 257$ &
-0.118800 & 2.0664 & (0.5335,0.5639) & \cite{sahin2003} \\
\hline
 &$\vNChz \times \Pcf$ & $256\times 256$ &
-0.122151 & 1.98912 & (0.5215,0.5449) & \eqref{cavityBD} \\
& $\vcQqhz \times \dcQlhz$ & $128\times 128$ &
-0.121492 & 1.97645 & (0.5195,0.5430) & \eqref{watertight} \\
\rb{2500} & $\vcQqhz \times \dcQlhz$ & $128\times 128$ &
-0.115717 & 1.88476 & (0.5195,0.5430) & \eqref{leaky} \\
& \cite{erturk2005numerical} & $601 \times 601$ &
-0.121035 & 1.96968 & (0.5200,0.5433) & - \\
\hline 
&$\vNChz \times \Pcf$ & $256\times 256$ &
-0.122713 & 1.97778 & (0.5176,0.5410) & \eqref{cavityBD} \\
& $\vcQqhz \times \dcQlhz$ & $128\times 128$ &
-0.121860 & 1.96186 & (0.5195,0.5391) & \eqref{watertight} \\
& $\vcQqhz \times \dcQlhz$ & $128\times 128$ &
-0.115310 & 1.85833 & (0.5195,0.5430) & \eqref{leaky} \\
\rb{3200} & \cite{ghia} & $257\times 257$ &
-0.120377 & 1.98860 & (0.5165,0.5469) & - \\
& \cite{hna09} & $128\times 128$ &
-0.121768 & - & (0.5165,0.5352) & \cite{hna09} \\
& \cite{sahin2003} & $257\times 257$ &
-0.121628 & 1.9593 & (0.5201,0.5376) & \cite{sahin2003} \\
\hline
 & $\vNChz \times \Pcf$ & $256\times 256$ &
-0.123658 & 1.96650 & (0.5137,0.5371) & \eqref{cavityBD} \\
& $\vcQqhz \times \dcQlhz$ & $128\times 128$ &
-0.122368 & 1.94277 & (0.5156,0.5352) & \eqref{watertight} \\
& $\vcQqhz \times \dcQlhz$ & $128\times 128$ &
-0.114120 & 1.81321 & (0.5156,0.5352) & \eqref{leaky} \\
5000 & \cite{erturk2005numerical} & $601 \times 601$ &
-0.121289 & 1.92660 & (0.5150,0.5350) & - \\
& \cite{ghia} & $257\times 257$ &
-0.118966 & 1.86016 & (0.5117,0.5352) & - \\
& \cite{hna09} & $128\times 128$ &
-0.121218 & - & (0.5156,0.5352) & \cite{hna09} \\
& \cite{sahin2003} & $257\times 257$ &
-0.122050 & 1.9392 & (0.5134,0.5376) & \cite{sahin2003} \\ \hline
\end{tabular}
\end{scriptsize}
\end{table}

\begin{table}
\caption{Computed secondary vortices;
the values of stream function ($\psi$) and location $(x,y)$
for the $\vNChz \times \Pcf$ element.}
\label{secondvortex}
\centering
\begin{scriptsize}
\begin{tabular}{ccccccc}
 & \multicolumn{2}{c}{Bottom left}
 & \multicolumn{2}{c}{Bottom Right}
 & \multicolumn{2}{c}{Top left} \\ \hline
Re & $\psi_{\max}$ & $(x,y)$ & $\psi_{\max}$ & $(x,y)$ &
$\psi_{\max}$ & $(x,y)$ \\ \hline
100 & 1.7368E-06 & (0.0332,0.0332)
& 1.2597E-05 & (0.9434,0.0605)
& - & - \\
400 & 1.4100E-05 & (0.0488,0.0488)
& 6.4495E-04 & (0.8848,0.1230)
& - & - \\
1000 & 2.3223E-04 & (0.0840,0.0762) 
& 1.7319E-03 & (0.8652,0.1113)
& - & - \\
2500 & 9.2779E-04 & (0.0840,0.1113)
& 2.6661E-03 & (0.8340,0.0918)
& 3.3918E-04 & (0.0410,0.8887) \\
3200 & 1.1104E-03 & (0.0801,0.1191)
& 2.8323E-03 & (0.8223,0.0840)
& 7.0750E-04 & (0.0527,0.8965) \\
5000 & 1.3660E-03 & (0.0723,0.1387)
& 3.0641E-03 & (0.8027,0.0723)
& 1.4566e-03 & (0.0645,0.9082) \\ \hline
\end{tabular}
\end{scriptsize}
\end{table}

In \secref{sec:acc}, we introduced the indicators for the accuracy
of the numerical solution.
First, the volumetric flow rate values $Q_{u,x_{c}}$ and $Q_{v,y_{c}}$ 
defined by \eqref{eq:volumetric} are shown in \tabref{volflowrate}.
for the $\vNChz \times \Pcf$ element 
The values of $Q_{u,x_{c}}$ and $Q_{v,y_{c}}$  for the $\vcQqhz \times \dcQlhz$
element are much larger than those values at $Q_{u,x_{c}-h/2}$,
$Q_{u,x_{c}+h/2}$, $Q_{v,y_{c}-h/2}$, and $Q_{v,y_{c}+h/2}$
for the $\vNChz \times \Pcf$ element.
Erturk {\it et al.} \cite{erturk2005numerical} calculated $Q_{u,x_{c}}$and $Q_{v,y_{c}}$ by using
their solutions. The smallest values of $Q_{u,x_{c}}=4.5\text{E-8}$ and
$Q_{v,y_{c}}=1.34\text{E-7}$ in \cite{erturk2005numerical} are lager than the largest values of
$Q_{u}$ and $Q_{v}$ for the $\vNChz \times \Pcf$ element.

\tabref{iaccuracy} shows the values of \eqref{vorcondition} and \eqref{qdiv}
for the $\vNChz \times \Pcf$ and the $\vcQqhz \times \dcQlhz$ elements.
Concerning the compatibility condition \eqref{vorcondition}, the $\vNChz
\times \Pcf$
element and the Taylor-Hood element with the leaky cavity boundary condition
\eqref{leaky} give precise values, while 
the Taylor-Hood element with the watertight cavity boundary condition
\eqref{watertight} generates about 0.3\% errors.
An investigation of \eqref{qdiv} shows that the numerical results obtained by
using the $\vNChz \times \Pcf$ element
are more accurate than those by the Taylor-Hood element.
Moreover, the absolute values \eqref{qdiv} for the $\vNChz \times \Pcf$ element are
independent of Reynolds number and element $Q_{jk}$ due to \thmref{thm:div-q}.
With the grid size $256 \times 256$, the absolute values \eqref{qdiv} for the 
$\vNChz \times \Pcf$ element is given by
\beq
\max_{Q_{jk} \in \Tau_h} \left|
\int_{Q_{jk}} \div \bu_h \dx
\right|
=\frac{1}{256^{3}} \approx 5.9605\text{E-8},
\eeq
while such values for the Taylor-Hood element with watertight and
leaky cavity boundary conditions are given in \tabref{iaccuracy}.
It should be stressed that the values obtained by the $\vNChz \times \Pcf$
element are smaller by a factor of four than those obtained by the Taylor-Hood element.

At least judged by the three accuracy indicators, \eqref{eq:volumetric},
\eqref{vorcondition}, and \eqref{qdiv}, the numerical solutions
by using the $\vNChz \times \Pcf$ element without any modification at the top
corners are more accurate than those by using the
Taylor-Hood element with modified boundary conditions \eqref{watertight} and \eqref{leaky}.

\begin{table}
\caption{Volumetric flow rates along the vertical and horizontal lines through
the geometric center of the cavity, $(x_{c},y_{c})$, by using the $\vNChz \times \Pcf$ and
$\vcQqhz \times \dcQlhz$ element}
\label{volflowrate}
\centering
\begin{tabular}{c|ccccc}
\hline
Re & FEM & Grid & $Q_{u,x_{c}}$ & $Q_{v,y_{c}}$  & BC \\ \hline
 & $\vNChz \times \Pcf$ & $256\times 256$ &
1.9039e-16 & 1.2514e-13 & \eqref{cavityBD} \\
\rb{100} & $\vcQqhz \times \dcQlhz$ & $128\times 128$ &
9.3009E-06 & 6.5662E-08 & \eqref{watertight} \\
 & $\vcQqhz \times \dcQlhz$ & $128\times 128$ &
1.3114E-03 & 9.7804E-08 & \eqref{leaky} \\
\hline
 & $\vNChz \times \Pcf$ & $256\times 256$ &
2.1554e-16 & 1.3347e-13 & \eqref{cavityBD} \\
\rb{400} & $\vcQqhz \times \dcQlhz$ & $128\times 128$ &
1.4876E-05 & 1.2495E-06 & \eqref{watertight} \\
 & $\vcQqhz \times \dcQlhz$ & $128\times 128$ &
1.3170E-03 & 1.1132E-06 & \eqref{leaky} \\
\hline
 & $\vNChz \times \Pcf$ & $256\times 256$ &
3.5996e-17 & 1.1037e-14 & \eqref{cavityBD} \\
\rb{1000} & $\vcQqhz \times \dcQlhz$ & $128\times 128$ &
2.4097E-05 & 2.8794E-06 & \eqref{watertight} \\
 & $\vcQqhz \times \dcQlhz$ & $128\times 128$ &
1.3264E-03 & 2.5407E-06 & \eqref{leaky} \\
\hline
 & $\vNChz \times \Pcf$ & $256\times 256$ &
2.4373e-16 & 1.5280e-13 & \eqref{cavityBD} \\
\rb{2500} & $\vcQqhz \times \dcQlhz$ & $128\times 128$ &
4.0694E-05 & 5.7553E-06 & \eqref{watertight} \\
 & $\vcQqhz \times \dcQlhz$ & $128\times 128$ &
1.3431E-03 & 4.8544E-06 & \eqref{leaky} \\
\hline
 & $\vNChz \times \Pcf$ & $256\times 256$ &
2.1814e-16 & 5.1092e-14 & \eqref{cavityBD} \\
\rb{3200} & $\vcQqhz \times \dcQlhz$ & $128\times 128$ &
4.6986E-05 & 7.0223E-06 & \eqref{watertight} \\
 & $\vcQqhz \times \dcQlhz$ & $128\times 128$ &
1.3494E-03 & 5.8405E-06 & \eqref{leaky} \\
\hline
 & $\vNChz \times \Pcf$ & $256\times 256$ &
3.5562e-16 & 1.2311e-13 & \eqref{cavityBD} \\ 
\rb{5000} & $\vcQqhz \times \dcQlhz$ & $128\times 128$ &
6.1055E-05 & 1.0206E-05 & \eqref{watertight} \\
 & $\vcQqhz \times \dcQlhz$ & $128\times 128$ &
1.3634E-03 & 8.2691E-06 & \eqref{leaky} \\ \hline
\end{tabular}
\end{table}

\begin{table}
\caption{Compatibility \eqref{vorcondition}  and
  incompressibility conditions \eqref{qdiv} 
for the $\vNChz \times \Pcf$ and $\vcQqhz \times \dcQlhz$ elements.} \label{iaccuracy}
\centering
\begin{tabular}{c|ccccc}
\hline
Re & FEM & Grid & $| \int_{\O} \omega  \; \dx + 1 |$ &\eqref{qdiv} & BC \\ \hline
 & $\vNChz \times \Pcf$ & $256 \times 256$ & 2.8866e-15 & 5.9605E-08 & \eqref{cavityBD} \\
100 & $\vcQqhz \times \dcQlhz$ & $128\times 128$ & 2.6042e-03 & 6.1596E-04 & \eqref{watertight} \\
& $\vcQqhz \times \dcQlhz$ & $128\times 128$ & 1.1102e-15 & 3.3407E-04 &
\eqref{leaky} \\
\hline
 & $\vNChz \times \Pcf$ & $256 \times 256$ & 2.2204e-16 & 5.9605E-08 & \eqref{cavityBD} \\
400 & $\vcQqhz \times \dcQlhz$ & $128\times 128$ & 2.6042e-03 & 6.6730E-04 & \eqref{watertight} \\
& $\vcQqhz \times \dcQlhz$ & $128\times 128$ & 4.7740e-15 & 3.9730E-04 &
\eqref{leaky} \\
\hline 
 & $\vNChz \times \Pcf$ & $256 \times 256$ & 2.6645e-15 & 5.9605E-08 & \eqref{cavityBD} \\
1000 & $\vcQqhz \times \dcQlhz$ & $128\times 128$ & 2.6042e-03 & 7.2274E-04 & \eqref{watertight} \\
& $\vcQqhz \times \dcQlhz$ & $128\times 128$ & 1.5543e-15 & 5.0746E-04 &
\eqref{leaky} \\
\hline 
 & $\vNChz \times \Pcf$ & $256 \times 256$ & 1.1102e-15 & 5.9605E-08 & \eqref{cavityBD} \\
2500 & $\vcQqhz \times \dcQlhz$ & $128\times 128$ & 2.6042e-03 & 1.1836E-03 & \eqref{watertight} \\
& $\vcQqhz \times \dcQlhz$ & $128\times 128$ & 7.7716e-15 & 5.9441E-04 &
\eqref{leaky} \\
\hline 
 & $\vNChz \times \Pcf$ & $256 \times 256$ & 3.9968e-15 & 5.9605E-08 & \eqref{cavityBD} \\
3200 & $\vcQqhz \times \dcQlhz$ & $128\times 128$ & 2.6042e-03 & 1.3685E-03 & \eqref{watertight} \\
& $\vcQqhz \times \dcQlhz$ & $128\times 128$ & 1.5543e-15 & 6.0657E-04 &
\eqref{leaky} \\
\hline 
 & $\vNChz \times \Pcf$ & $256 \times 256$ & 1.4433e-15 & 5.9605E-08 & \eqref{cavityBD} \\
5000 & $\vcQqhz \times \dcQlhz$ & $128\times 128$ & 2.6042e-03 & 1.6240E-03 & \eqref{watertight} \\
& $\vcQqhz \times \dcQlhz$ & $128\times 128$ & 4.2188e-15 & 6.7909E-04 & \eqref{leaky} \\ \hline
\end{tabular}
\end{table}

\section{Conclusions} \label{sec:conc}
%

The $\vNChz \times \Pcf$ element is applied to
solve the lid driven cavity problem with least modification at the two top
corner element to deal with the jump discontinuities there using 
the DSSY element (of CDY element).

The numerical solutions using $\vNChz \times \Pcf$ element are compared with
bench mark solutions and the horizontal and vertical components of the
velocity at the center are correct up to mostly two and three digits if the mesh
sizes are  $256\times256$ and $512\times 512,$ respectively.

Numerical solutions were compared with those
the conforming $\vcQqhz \times \dcQlhz$ element (Taylor-Hood element) with leaky and watertight
cavity boundary conditions.
Three indicators for accuracy of the numerical solution have been
compared.
(1) The incompressibility condition
(2) The compatibility condition 
(3) 
with the Neumann boundary condition are used to check the
accuracy of the numerical solutions.

Our numerical solutions satisfy the incompressibility and compatibility condition
precisely. Numerical results computed by using the $\vNChz \times \Pcf$ element
show the best results in terms of satisfying incompressibility and compatibility conditions,
and volumetric flow rates.

\section*{Acknowledgments}
The authors are very grateful to Prof. Roland Glowinski who inspired us to
investigate in this approach to treat the corner singularities in the
approximation of lid cavity flows. Also, the work has been initiated while the
second author was visiting Texas A\&M
University. He thanks the Department of Mathematics and the Institute for Scientific
Computation of Texas A\&M University for financial and other administrative  supports during his visit.


\begin{thebibliography}{10}

\bibitem{ifiss}
Incompressible {F}low \& {I}terative {S}olver {S}oftware.
\newblock {http://www.maths.}{manchester.ac.uk}{/~djs/ifiss}.

\bibitem{altmann-carstensen}
R.~Altmann and C.~Carstensen.
\newblock {$P_1$}-nonconforming finite elements on triangulations into
  triangles and quadrilaterals.
\newblock {\em SIAM J. Numer. Anal.}, 50(2):418--438, 2011.

\bibitem{altmann2012p}
R.~Altmann and C.~Carstensen.
\newblock {$P_1$}-nonconforming finite elements on triangulations into triangles
  and quadrilaterals.
\newblock {\em SIAM Journal on Numerical Analysis}, 50(2):418--438, 2012.

\bibitem{auteri2002numerical}
F.~Auteri, N.~Parolini, and L.~Quartapelle.
\newblock Numerical investigation on the stability of singular driven cavity
  flow.
\newblock {\em Journal of Computational Physics}, 183(1):1--25, 2002.

\bibitem{aydin}
M.~Aydin and R.~Fenner.
\newblock Boundary element analysis of driven cavity flow for low and moderate
  {Reynolds} number.
\newblock {\em Int. J. Numer. Meth. Fluids.}, 37:45--64, 2001.

\bibitem{barragy1997}
E.~Barragy and G.~Carey.
\newblock Stream function-vorticity driven cavity solution using $p$ finite
  elements.
\newblock {\em Computers \& Fluids}, 26:453--468, 1997.

\bibitem{bercovier-pironneau-79}
M.~Bercovier and O.~Pironneau.
\newblock Error estimates for finite element method solution of the {Stokes}
  problem in the primitive variables.
\newblock {\em Numer. Math.}, 33(2):211--224, 1979.

\bibitem{botella1998benchmark}
O.~Botella and R.~Peyret.
\newblock Benchmark spectral results on the lid-driven cavity flow.
\newblock {\em Computers \& Fluids}, 27(4):421--433, 1998.

\bibitem{bruneau20062d}
C.-H. Bruneau and M.~Saad.
\newblock The {2D} lid--driven cavity problem revisited.
\newblock {\em Computers \& Fluids}, 35(3):326--348, 2006.

\bibitem{cdssy}
Z.~Cai, J.~{Douglas,~Jr.}, J.~E. Santos, D.~Sheen, and X.~Ye.
\newblock Nonconforming quadrilateral finite elements: {A} correction.
\newblock {\em Calcolo}, 37(4):253--254, 2000.

\bibitem{cdy}
Z.~Cai, J.~{Douglas,~Jr.}, and X.~Ye.
\newblock A stable nonconforming quadrilateral finite element method for the
  stationary {S}tokes and {N}avier-{S}tokes equations.
\newblock {\em Calcolo}, 36:215--232, 1999.

\bibitem{cai-wang-cavity}
Z.~Cai and Y.~Wang.
\newblock An error estimate for two-dimensional {Stokes} driven cavity flow.
\newblock {\em Math. Comp.}, 78:771--787, 2008.

\bibitem{crouzeix-raviart}
M.~Crouzeix and P.-A. Raviart.
\newblock Conforming and nonconforming finite element methods for solving the
  stationary {Stokes} equations.
\newblock {\em R.A.I.R.O.-- Math. Model. Anal. Numer.}, 7:33--75, 1973.

\bibitem{cuvelier1986finite}
C.~Cuvelier, A.~Segal, and A.~A. Van~Steenhoven.
\newblock {\em Finite element methods and {Navier}--{Stokes} equations},
  volume~22.
\newblock Springer, 1986.

\bibitem{dssy-nc-ell}
J.~{Douglas,~Jr.}, J.~E. Santos, D.~Sheen, and X.~Ye.
\newblock Nonconforming {G}alerkin methods based on quadrilateral elements for
  second order elliptic problems.
\newblock {\em ESAIM--Math. Model. Numer. Anal.}, 33(4):747--770, 1999.

\bibitem{elman2014finite}
H.~Elman, D.~Silvester, and A.~Wathen.
\newblock {\em Finite elements and fast iterative solvers: with applications in
  incompressible fluid dynamics}.
\newblock Oxford University Press, 2014.

\bibitem{erturk2009discussions}
E.~Erturk.
\newblock Discussions on driven cavity flow.
\newblock {\em International Journal for Numerical Methods in Fluids},
  60(3):275--294, 2009.

\bibitem{erturk2005numerical}
E.~Erturk, T.~C. Corke, and C.~G{\"o}k{\c{c}}{\"o}l.
\newblock Numerical solutions of {2-D} steady incompressible driven cavity flow
  at high {Reynolds} numbers.
\newblock {\em International Journal for Numerical Methods in Fluids},
  48(7):747--774, 2005.

\bibitem{ghia}
U.~Ghia, K.~N. Ghia, and C.~T. Shin.
\newblock High-{Re} solutions for incompressible flow using the
  {Navier}-{Stokes} equations and a multigrid method.
\newblock {\em J. Comp. Phys.}, 48:387--411, 1982.

\bibitem{hna09}
R.~Glowinski.
\newblock Finite element methods for incompressible viscous flow.
\newblock In P.~G. Ciarlet and J.~L. Lions, editors, {\em Handbook of Numerical
  Analysis. IX. Numerical Methods for Fluids (Part 3)}. Elsevier/North-Holland,
  Amsterdam, 2003.

\bibitem{glowinski-semicircular-06}
R.~Glowinski, G.~Guidoboni, and T.-W. Pan.
\newblock Wall-driven incompressible viscous flow in a two-dimensional
  semi-circular cavity.
\newblock {\em J. Comp. Phys.}, 216(1):76--91, 2006.

\bibitem{guermond2011new}
J.-L. Guermond and P.~Minev.
\newblock A new class of massively parallel direction splitting for the
  incompressible {Navier}--{Stokes} equations.
\newblock {\em Computer Methods in Applied Mechanics and Engineering},
  200(23):2083--2093, 2011.

\bibitem{guermond2010new}
J.-L. Guermond and P.~D. Minev.
\newblock A new class of fractional step techniques for the incompressible
  {Navier}--{Stokes} equations using direction splitting.
\newblock {\em Comptes Rendus Mathematique}, 348(9):581--585, 2010.

\bibitem{guermond2012start}
J.-L. Guermond and P.~D. Minev.
\newblock Start-up flow in a three-dimensional lid-driven cavity by means of a
  massively parallel direction splitting algorithm.
\newblock {\em International Journal for Numerical Methods in Fluids},
  68(7):856--871, 2012.

\bibitem{taylor-hood}
P.~Hood and C.~Taylor.
\newblock A numerical solution of the {Navier}--{Stokes} equations using the
  finite element techniques.
\newblock {\em Computers \& Fluids}, 1:73--100, 1973.

\bibitem{jeon-nam-sheen-shim-nonpara}
Y.~Jeon, H.~Nam, D.~Sheen, and K.~Shim.
\newblock A class of nonparametric {DSSY} nonconforming quadrilateral elements.
\newblock {\em ESAIM--Math. Model. Numer. Anal.}, 47(06):1783--1796, 2013.

\bibitem{karakashian1982galerkin}
O.~A. Karakashian.
\newblock On a {Galerkin}--{Lagrange} multiplier method for the stationary
  {Navier}--{Stokes} equations.
\newblock {\em SIAM J. Numer. Anal.}, 19(5):909--923, 1982.

\bibitem{stab-cheapest}
S.~Kim, J.~Yim, and D.~Sheen.
\newblock Stable cheapest nonconforming finite elements for the {Stokes}
  equations.
\newblock {\em J. Comput. Appl. Math.}, 299:2-14, 2016.

\bibitem{cpark-thesis}
C.~Park.
\newblock {\em A study on locking phenomena in finite element methods}.
\newblock PhD thesis, Department of Mathematics, Seoul National University,
  Korea, Feb. 2002.
\newblock Available at http://www.nasc.snu.ac.kr/cpark/papers/phdthesis.ps.gz.

\bibitem{parksheen-p1quad}
C.~Park and D.~Sheen.
\newblock {$P_1$}-nonconforming quadrilateral finite element methods for
  second-order elliptic problems.
\newblock {\em SIAM J. Numer. Anal.}, 41(2):624--640, 2003.

\bibitem{rann}
R.~Rannacher and S.~Turek.
\newblock Simple nonconforming quadrilateral {Stokes} element.
\newblock {\em Numer. Methods Partial Differential Equations}, 8:97--111, 1992.

\bibitem{sahin2003}
M.~Sahin and R.~Owens.
\newblock A novel fully implicit finite volume methods applied to the
  lid-driven cavity problem--{Part I:} {High} {Reynolds} number flow
  calculations.
\newblock {\em Int. J. Numer. Meth. Fluids.}, 42:57--77, 2003.

\bibitem{shen1991}
J.~Shen.
\newblock Hopf bifurcation of the unsteady regularized driven cavity flow.
\newblock {\em J. Comp. Phys.}, 95:228--245, 1991.

\end{thebibliography}

\end{document}